\documentclass[preprint,12pt]{elsarticle}
\usepackage{amsfonts}
\usepackage{graphics}
\usepackage{graphicx}
\usepackage{amsmath,amssymb}
\usepackage{bm}
\usepackage{soul}
\usepackage[normalem]{ulem}
\usepackage{float}

\usepackage{enumitem}

\biboptions{sort&compress}

\usepackage{tikz}
\usetikzlibrary{calc,decorations.pathreplacing}
\usetikzlibrary{decorations.markings}

\numberwithin{equation}{section}

\allowdisplaybreaks

\journal{Indagationes Mathematicae}

\begin{document}

\begin{frontmatter}

\title{The realm of finite lattices in combination with a new dimension}

\author[label1]{D. N. Georgiou\corref{cor1}}\ead{georgiou@math.upatras.gr}
\address[label1]{University of Patras, Department of Mathematics, 26504 Patra, Greece}
\cortext[cor1]{Corresponding Author}

\author[label2]{Y. Hattori}\ead{hyasu2014@gmail.com}
\address[label2]{Shimane University (Professor Emeritus), Matsue, 690-8504, Japan}

\author[label1]{A. C. Megaritis}\ead{acmegaritis@upatras.gr}

\author[label1]{F. Sereti}\ead{seretifot@gmail.com}

\begin{abstract}
The Ordered Set Theory is a branch of Mathematics that studies
partially ordered sets (usually posets) and lattices. The meaning of
dimension is one of the main parts of this field. Dimensions of
partially ordered sets and lattices have been studied in various
researches. In particular, the covering dimension, the Krull
dimension and the small inductive dimension have been studied
extensively for the class of finite lattices. In this paper, we
insert a new meaning of dimension for finite lattices called large
inductive dimension and denoted by Ind. We study various of its
properties based on minimal covers. Also, given two finite lattices,
we study the dimension Ind of their linear sum, Cartesian,
lexicographic and rectangular product, investigating the ``behavior"
of this dimension. In addition, we study relations of this new
dimension with the small inductive dimension, covering dimension and
Krull dimension, presenting various facts and examples that
strengthen the corresponding results.
\end{abstract}


\end{frontmatter}

\newtheorem{theorem}{Theorem}[section]
\newtheorem{proposition}[theorem]{Proposition}
\newtheorem{lemma}[theorem]{Lemma}
\newtheorem{corollary}[theorem]{Corollary}
\newtheorem{example}[theorem]{Example}
\newtheorem{remark}[theorem]{Remark}
\newtheorem{definition}[theorem]{Definition}
\newtheorem{notation}[theorem]{Notation}
\newtheorem{question}[theorem]{Question}

\section{Introduction}\label{sec1}

Partially ordered sets (in short posets), lattices and frames are
the main part of many researches. They form a constantly growing
chapter which attracts the interest of many research teams (see for
example \cite{DMa1, Du, DMNa, WTM, WJi}), investigating new
properties of posets. Simultaneously, there is no doubt that
dimensions of topological spaces are the base of many studies. There
are many results for the covering dimension, the small inductive
dimension and the large inductive dimension in various classes of
topological spaces. We refer to \cite{BW,C2,E,N1,N2,P} for more
studies on these dimensions.

Recently, the covering dimension, the small inductive dimension and
the large inductive dimension had a particular interest in the class
of frames (see for example \cite{BG,BB1,BB2,C,C1}). In addition,
several researches are devoted to the notion of dimension for finite
partially ordered sets and finite lattices (see for example
\cite{Bae1,DP,DGMM,DGMS,GMPS,GMS,T1,T2,V,YM,ZZZ,WWY}). The order
dimension, the Krull dimension, the covering dimension and the small
inductive dimension developed their own significant chapter in
Dimension Lattice Theory (see for example \cite{BGMS1,GMPS,T2,ZZZ}).

In this paper, we study the notion of the large inductive dimension,
$\mathrm{Ind}$, in the class of finite lattices, proving properties
of this dimension, like the sublattice, sum and product properties.
The motivation of this study is based on the fact that posets and
lattices are used as models to investigate properties of digital
images (see for example \cite{EKM,KKM}) and generally, topological
properties of digital images and discrete arrays in two or more
dimensions. Thus, the notion of dimension has its special study and
a research issue that arises is the study of dimensions for finite
lattices.

The paper is organized as follows. In Section \ref{sec2}, we present
definitions that are used in the following sections. In Section
\ref{sec3}, we present the meaning of the large inductive dimension
for the class of finite lattices and some basic results. In Section
\ref{sec4}, we compare the large inductive dimension with the small
inductive dimension, proving that the large inductive dimension is
always greater than or equal to the small inductive dimension. Also,
we study the ``gap" between these two dimensions. In addition, we
present facts and various examples, investigating the relation
between the large inductive dimension and the covering dimension as
well as the Krull dimension and the height. Finally, in Section
\ref{sec5}, given two finite lattices, we study the large inductive
dimension of their linear sum, Cartesian product, lexicographic
product and rectangular product.

\section{Preliminaries}\label{sec2}

In this section we recall basic definitions and notations which will
be used in order to present our study. Throughout this paper we
suppose that the lattices are finite.

Let $L$ be a lattice. We denote by $0_L$ and $1_L$ the bottom and
the top element of $L$, respectively. For an element $x\in L$, we
use the notation ${\uparrow}x=\{y\in L:x\leqslant y\}$. Since the
dimension Ind is defined inductively for finite lattices (see
Definition \ref{Ind} below), we remind the reader that a
\emph{sublattice} of $L$ is a subset of $L$, which is a lattice with
the same meet and join operations as $L$. Every set ${\uparrow}x$ is
a sublattice of $L$.

Moreover, we use the notation $x^{*}=\bigvee\{y\in L:x\wedge
y=0_L\}$. In the case where $x^{*}\wedge x=0_L$, that is
$x^{*}=\max\{y\in L:y\wedge x=0_L\}$, the element $x^{*}$ is called
\emph{pseudocomplement} of $x$.

It is known that every finite distributive lattice is
pseudocomplemented. Generally, the above notion of pseudocomplement
plays an important role in related studies with several properties
on such lattices. However, in this paper we will adopt only the
asterisk notation and the term ``pseudocomplement" in order to have
the same terminology as in \cite{BB1}. Since the lattices we deal
with in this paper are not generally distributive, the use of the
term ``pseudocomplement" will be an agreement for the following
sections, obtaining that it is an element taking as $\bigvee\{y\in
L:x\wedge y=0_L\}$.

A subset $V$ of a lattice $L$ is called a \emph{cover} of $L$ if
$0_L \notin V$ and $\bigvee V=1_L$. A cover $U$ of a lattice $L$ is
called a \emph{refinement} of a cover $V$ of $L$, writing $U\succ
V$, if for each $u\in U$, there exists $v\in V$ such that
$u\leqslant v$. Especially, a cover $V$ of $L$ is said to be a
\emph{minimal cover} of $L$ if $V\subseteq C$ for every refinement
$C$ of $V$. (We state that in related bibliography (see for example
\cite{FJN}) we can find variations of the meaning of minimal covers
such as the minimal $\ll$-representations of the top element of $L$
consisting of join-irreducible elements of $L$. However, in our
study we keep in mind the first definition of minimal covers.)

Finally, we remind the reader that a subset $S$ of a lattice $L$ is
called \emph{down-set} if for every $x\in L$ with $x\leqslant s$ for
some $s\in S$, then $x\in S$.

In the whole paper we denote by $\mathbb{N}$ the set of natural
numbers and whenever it is necessary we write $\{1,2,\ldots\}$ for
the set of natural numbers except zero.

\section{The large inductive dimension and finite lattices}\label{sec3}

In this section, we study the large inductive dimension
$\mathrm{Ind}$ in the class of finite lattices and give basic
properties of this dimension.

\begin{definition} {\rm \label{Ind} Let $L$ be a finite lattice. The \emph{large inductive dimension}, $\mathrm{Ind}$, of $L$ is defined as
follows:
\begin{enumerate}
\item $\mathrm{Ind}(L)=-1$ if and only if $L=\{0_L\}$.
\item $\mathrm{Ind}(L)\leqslant k$, where $k\in \mathbb{N}$, if
for every $a\in L$ and for every $v\in L$ such that $a\vee v=1_L$,
there exists $u\in L$ such that $u\leqslant v$, $a\vee u=1_L$ and
$\mathrm{Ind} ({\uparrow}(u^{*} \vee u))\leqslant k-1$.
\item $\mathrm{Ind}(L)=k$, where $k\in \mathbb{N}$, if $\mathrm{Ind}(L)\leqslant k$ and $\mathrm{Ind}(L)\nleqslant
k-1$.
\end{enumerate}}
\end{definition}

Analogously to Proposition 4.6 of \cite{BB1}, isomorphic lattices
have the same large inductive dimension. Also, the simple lattice
$L=\{0_L,1_L\}$ has trivially $\mathrm{Ind}(L)=0$.

\begin{example}{\rm \label{Ex1} We consider the finite lattices represented by the
diagrams of Figure \ref{fig:1}.
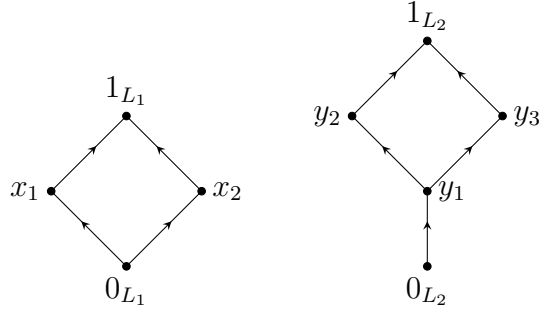
\begin{figure}[H]
\tikzstyle arrowstyle=[scale=1] \tikzstyle
directed=[postaction={decorate,decoration={markings, mark=at
position .6 with {\arrow[arrowstyle]{stealth}}}}]
\begin{center}
\begin{tikzpicture}
\draw[directed] (0,0)--(-1,1); \draw[directed] (0,0)--(1,1);
\draw[directed] (-1,)--(0,2); \draw[directed] (1,1)--(0,2);
\draw[fill] (0,0) circle [radius=0.05] node [below] {$0_{L_{1}}$};
\draw[fill] (-1,1) circle [radius=0.05] node [left] {$x_1$};
\draw[fill] (1,1) circle [radius=0.05] node [right] {$x_2$};
\draw[fill] (0,2) circle [radius=0.05] node [above] {$1_{L_{1}}$};
\draw[directed] (4,0)--(4,1); \draw[directed] (4,1)--(3,2);
\draw[directed] (4,1)--(5,2); \draw[directed] (3,2)--(4,3);
\draw[directed] (5,2)--(4,3); \draw[fill] (4,0) circle [radius=0.05]
node [below] {$0_{L_{2}}$}; \draw[fill] (4,1) circle [radius=0.05]
node [right] {$y_1$}; \draw[fill] (3,2) circle [radius=0.05] node
[left] {$y_2$}; \draw[fill] (5,2) circle [radius=0.05] node [right]
{$y_3$}; \draw[fill] (4,3) circle [radius=0.05] node [above]
{$1_{L_{2}}$};
\end{tikzpicture}
\end{center}
\caption{The finite lattices $(L_1,\leqslant_1)$ and $(L_2,\leqslant
_2)$} \label{fig:1}
\end{figure}
We have that $\mathrm{Ind}(L_1)=0$ and $\mathrm{Ind} (L_2)=1$. For
example, we apply Definition \ref{Ind} to see that $\mathrm{Ind}
(L_2)=1$.

(1) For $a=0_{L_2}$ and for every $v\in L_2$ with $a\vee
v=1_{L_{2}}$, that is $v=1_{L_{2}}$, there exists $u=v=1_{L_{2}}$
with $a\vee u=1_{L_{2}}$ and $${\rm
Ind}({\uparrow}((1_{L_{2}})^{*}\vee 1_{L_{2}}))={\rm
Ind}({\uparrow}(0_{L_{2}}\vee 1_{L_{2}}))={\rm
Ind}({\uparrow}1_{L_{2}})={\rm Ind}(\{1_{L_{2}}\})=-1.$$

(2) For $a=y_1$ and for every $v\in L_2$ with $a\vee v=1_{L_{2}}$,
that is $v=1_{L_{2}}$, there exists $u=v=1_{L_{2}}$ with $a\vee
u=1_{L_{2}}$ and $${\rm Ind}({\uparrow}((1_{L_{2}})^{*}\vee
1_{L_{2}}))={\rm Ind}({\uparrow}(0_{L_{2}}\vee 1_{L_{2}}))={\rm
Ind}({\uparrow}1_{L_{2}})={\rm Ind}(\{1_{L_{2}}\})=-1.$$

(3) For $a=y_2$ and for every $v\in L_2$ with $a\vee v=1_{L_{2}}$,
that is $v\in\{y_3,1_{L_{2}}\}$, there exists $u=y_3$ with
$u\leqslant v$, $a\vee u=1_{L_{2}}$ and $${\rm
Ind}({\uparrow}(y_3^{*}\vee y_3))={\rm Ind}({\uparrow}(0_{L_{3}}\vee
y_3))={\rm Ind}({\uparrow}y_3)=0.$$

(4) For $a=y_3$ and for every $v\in L_2$ with $a\vee v=1_{L_{2}}$,
that is $v\in\{y_2,1_{L_{2}}\}$, there exists $u=y_2$ with
$u\leqslant v$, $a\vee u=1_{L_{2}}$ and
$${\rm Ind}({\uparrow}(y_2^{*}\vee y_2))={\rm
Ind}({\uparrow}(0_{L_{2}}\vee y_2))={\rm Ind}({\uparrow}y_2)=0.$$

(5) For $a=1_{L_{2}}$ and for every $v\in L_2$ with $a\vee
v=1_{L_{2}}$, that is $v$ can be any of the elements of $L_2$, there
exists $u=0_{L_{2}}$ with $u\leqslant v$, $a\vee u=1_{L_{2}}$ and
$${\rm Ind}({\uparrow}((0_{L_{2}})^{*}\vee 0_{L_{2}}))={\rm
Ind}({\uparrow}(1_{L_{2}}\vee 0_{L_{2}}))={\rm
Ind}({\uparrow}1_{L_{2}})={\rm Ind}(\{1_{L_{2}}\})=-1.$$}
\end{example}

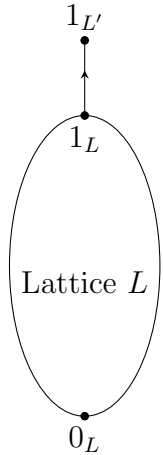
\begin{figure}[H]
\tikzstyle arrowstyle=[scale=1] \tikzstyle
directed=[postaction={decorate,decoration={markings,mark=at position
.6 with {\arrow[arrowstyle]{stealth}}}}]
\begin{center}
\begin{tikzpicture}

\draw[directed] (-1,4)--(-1,5);

\draw (-1,2) ellipse (1cm and 2cm);

\draw[fill] (-1,5) circle [radius=0.05] node [above]
{$1_{L^{\prime}}$}; \draw[fill] (-1,0) circle [radius=0.05] node
[below] {$0_L$}; \draw[fill] (-1,4) circle [radius=0.05] node
[below] {$1_L$}; \draw (-1,1.5) node [above] {Lattice $L$};
\end{tikzpicture}
\end{center}
\vspace{-0.5cm} \caption{The finite lattice $L^{\prime}$}
\label{FIG:2000}
\end{figure}

\begin{proposition} Let $L$ be a finite lattice and $k\in\mathbb{N}$ with ${\rm Ind}(L)=k$. Then there exists a finite lattice
$L^{\prime}$ with ${\rm Ind}(L^{\prime})=0$.
\end{proposition}
{\bf Proof.} We consider the finite lattice
$L^{\prime}=L\cup\{1_{L^{\prime}}\}$ of Figure \ref{FIG:2000}, where
$1_{L}\leqslant 1_{L^{\prime}}$, that is $1_{L^{\prime}}$ is the top
element of $L^{\prime}$. Then ${\rm Ind}(L^{\prime})=0$. Indeed, for
every $a\in L^{\prime}$ and for every $v\in L^{\prime}$ with $a\vee
v=1_{L^{\prime}}$, that is $v=1_{L^{\prime}}$, there exists
$u=v=1_{L^{\prime}}$ with $a\vee u=1_{L^{\prime}}$ and $${\rm
Ind}({\uparrow}((1_{L^{\prime}})^{*}\vee 1_{L^{\prime}}))={\rm
Ind}({\uparrow}(0_{L^{\prime}}\vee 1_{L^{\prime}}))={\rm
Ind}({\uparrow}1_{L^{\prime}})={\rm Ind}(\{1_{L^{\prime}}\})=-1.\
\Box$$

A modification of Definition \ref{Ind} (2) is given in the following
proposition.

\begin{theorem} Let $L$ be a finite lattice and
$k\in\mathbb{N}$. Then ${\rm Ind}(L)\leqslant k$ if and only if for
every $a\in L$ and every $u\in\min(\{x\in L:a\vee x=1_L\})$ we have
${\rm Ind}({\uparrow}(u^{*}\vee u))\leqslant k-1$.
\end{theorem}
{\bf Proof.} Let ${\rm Ind}(L)\leqslant k$, $a\in L$ and
$u\in\min(\{x\in L:a\vee x=1_L\})$. We shall prove that ${\rm
Ind}({\uparrow}(u^{*}\vee u))\leqslant k-1$. Clearly, $a\vee u=1_L$.
Thus, since ${\rm Ind}(L)\leqslant k$, there exists $w\in L$ such
that $w\leqslant u$, $w\vee a=1_L$ and ${\rm
Ind}({\uparrow}(w^{*}\vee w))\leqslant k-1$. Then $u=w$ and the
result is proved.

Conversely, we shall prove that ${\rm Ind}(L)\leqslant k$. Let $a\in
L$ and $v\in L$ with $a\vee v=1_L$. Then by assumption, every
$u\in\min(\{x\in L:a\vee x=1_L\})$ satisfying $u\leqslant v$ has the
properties of $a\vee u=1_L$ and ${\rm Ind}({\uparrow}(u^{*}\vee
u))\leqslant k-1$. (The element $u$ always exists since in the worst
case we can take $u=v$.) Hence, ${\rm Ind}(L)\leqslant k$. $\Box$

\begin{theorem}\label{mth} For any $k\in \{1,2,\ldots\}$, there exists
a finite lattice $L$ with $\mathrm{Ind}(L)=k$. Also, the
pseudocomplement of every element of $L\setminus\{0_{L}\}$ is
$0_{L}$.
\end{theorem}
{\bf Proof.} We prove the theorem by induction on $k$. Firstly, we
shall prove that there exists a finite lattice $L$, in which the
pseudocomplement of every element of $L\setminus\{0_{L}\}$ is
$0_{L}$, with $\mathrm{Ind}(L)=1$. Indeed, if we consider the
lattice $(L_2,\leqslant_2)$ of Example \ref{Ex1}, we have the
desired result.

We suppose that there exists a finite lattice $M$, in which the
pseudocomplement of every element of $M\setminus\{0_{M}\}$ is
$0_{M}$, with $\mathrm{Ind}(M)=k-1$, where $k>1$. We shall prove
that there exists a finite lattice $L$, in which the
pseudocomplement of every element of $L\setminus\{0_{L}\}$ is
$0_{L}$, with $\mathrm{Ind}(L)=k$.

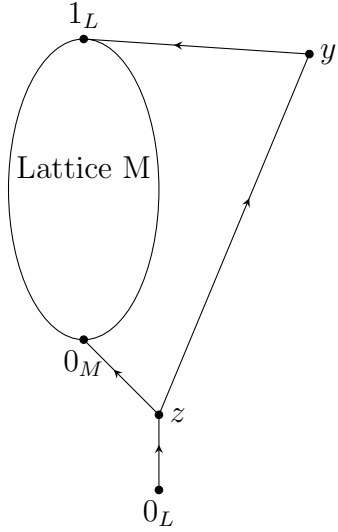
\begin{figure}[H]
\tikzstyle arrowstyle=[scale=1] \tikzstyle
directed=[postaction={decorate,decoration={markings,mark=at position
.6 with {\arrow[arrowstyle]{stealth}}}}]
\begin{center}
\begin{tikzpicture}
\draw[directed] (0,-3)--(0,-2); \draw[directed] (0,-2)--(-1,-1);
\draw[directed] (0,-2)--(2,2.8); \draw[directed] (2,2.8)--(-1,3);
\draw (-1,1) ellipse (1cm and 2cm); \draw[fill] (0,-3) circle
[radius=0.05] node [below] {$0_L$}; \draw[fill] (0,-2) circle
[radius=0.05] node [right] {$z$}; \draw[fill] (2,2.8) circle
[radius=0.05] node [right] {$y$}; \draw[fill] (-1,-1) circle
[radius=0.05] node [below] {$0_M$}; \draw[fill] (-1,3) circle
[radius=0.05] node [above] {$1_L$}; \draw (-1,1) node [above]
{Lattice M};
\end{tikzpicture}
\end{center}
\caption{The lattice $(L,\leqslant)$} \label{fig:2}
\end{figure}

We consider the finite lattice $L=\{0_L,z,y\} \cup M$ (see Figure
\ref{fig:2}). The partial order $\leqslant$ on $L$ is defined as
follows:
\begin{enumerate}
\item $z\leqslant w$, for every $w\in L \backslash \{0_L\}$,
\item $w\parallel y$, for every $w\in M$, that is,
$w\nleqslant y$ and $y\nleqslant w$.
\end{enumerate}

We observe that the pseudocomplement of every element of
$L\setminus\{0_L\}$ is $0_L$. We shall prove that
$\mathrm{Ind}(L)=k$. For that, we consider the following cases:

(i) For $a=0_L$ (respectively, $a=z$) we have $v=1_L$ with $a\vee
v=1_L$ and thus, there exists the element $u=1_L$ such that
$u\leqslant v$, $a\vee u=1_L$ and
\[\mathrm{Ind}({\uparrow}(u^{*} \vee u))=-1<k-1.\]

(ii) For $a=1_L$ and $v\in L$ with $a\vee v=1_L$, there exists the
element $u=0_L$ such that $u\leqslant v$, $a\vee u=1_L$ and
\[\mathrm{Ind}({\uparrow}(u^{*} \vee u))=-1<k-1.\]

(iii) For $a=0_M$ we have $v\in\{y,1_L\}$ with $a\vee v=1_L$ and
thus there exists the element $u=y$ such that $u\leqslant v$, $a\vee
u=1_L$ and \[\mathrm{Ind}({\uparrow}(y^{*} \vee
y))=\mathrm{Ind}({\uparrow}(0_L \vee
y))=\mathrm{Ind}({\uparrow}y)=0<k-1.\]

(iv) Let $a=w$, where $w\in M\setminus\{0_M,1_M\}$ (we state that
$1_M=1_L$) and $v\in M\cup\{y\}$ for which $a\vee v=1_L$. Then we
have the following cases:

(a) Let $a=w$ be any element of  $M\setminus\{0_M,1_M\}$ for which
$v\in\{y,1_L\}$. Then we consider the element $u=y$ such that
\[\mathrm{Ind}({\uparrow}(y^{*} \vee y))=\mathrm{Ind}({\uparrow}(0_L
\vee y))=\mathrm{Ind}({\uparrow}y)=0<k-1.\]

(b) Let $a=w$ be any element of  $M\setminus\{0_M,1_M\}$ for which
$v\in M$. Since $\mathrm{Ind}(M)=k-1$, there exists $u\in M$ such
that $u\leqslant v$, $a\vee u=1_L$ and
\[\mathrm{Ind}({\uparrow}(u^{*}_M \vee u))=\mathrm{Ind}({\uparrow}(0_M \vee
u))=\mathrm{Ind}({\uparrow}u)\leqslant k-2,\] and thus,
\[\mathrm{Ind}({\uparrow}(u^{*}_L \vee
u))=\mathrm{Ind}({\uparrow}(0_L \vee
u))=\mathrm{Ind}({\uparrow}u)\leqslant k-2.\]

(v) For $a=y$ we have $v\in M$ with $a\vee v=1_L$. Then there exists
$u=0_M$ such that $u\leqslant v$, $a\vee u=1_L$ and
\[\mathrm{Ind}({\uparrow}(0_M^{*} \vee 0_M))=\mathrm{Ind}({\uparrow}(0_L \vee
0_M))=\mathrm{Ind}({\uparrow}0_M)=\mathrm{Ind}(M)=k-1.\]

Therefore, $\mathrm{Ind}(L)\leqslant k$. Moreover, since
$\mathrm{Ind}(L)\nleqslant k-1$, we have that $\mathrm{Ind}(L)=k$.
$\Box$

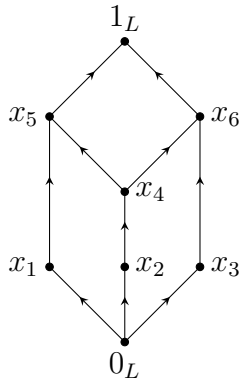
\begin{figure}[H]
\tikzstyle arrowstyle=[scale=1] \tikzstyle
directed=[postaction={decorate,decoration={markings,mark=at position
.6 with {\arrow[arrowstyle]{stealth}}}}]
\begin{center}
\begin{tikzpicture}
\draw[directed] (0,0)--(-1,1); \draw[directed] (0,0)--(0,1);
\draw[directed] (0,0)--(1,1); \draw[directed] (0,1)--(0,2);
\draw[directed] (-1,1)--(-1,3); \draw[directed] (1,1)--(1,3);
\draw[directed] (0,2)--(-1,3); \draw[directed] (0,2)--(1,3);
\draw[directed] (-1,3)--(0,4); \draw[directed] (1,3)--(0,4);
\draw[fill] (0,0) circle [radius=0.05] node [below] {$0_L$};
\draw[fill] (-1,1) circle [radius=0.05] node [left] {$x_1$};
\draw[fill] (0,1) circle [radius=0.05] node [right] {$x_2$};
\draw[fill] (1,1) circle [radius=0.05] node [right] {$x_3$};
\draw[fill] (0,2) circle [radius=0.05] node [right] {$x_4$};
\draw[fill] (-1,3) circle [radius=0.05] node [left] {$x_5$};
\draw[fill] (1,3) circle [radius=0.05] node [right] {$x_6$};
\draw[fill] (0,4) circle [radius=0.05] node [above] {$1_L$};
\end{tikzpicture}
\end{center}
\caption{The finite lattice $(L,\leqslant)$} \label{fig:3}
\end{figure}

\medskip
We complete this section, presenting a sublattice property for the
large inductive dimension. As we observe in Figure \ref{fig:3}, a
corresponding sublattice theorem does not hold without mentioning
further conditions for the sublattice of a finite lattice.

Especially, for the lattice of Figure \ref{fig:3} we have that
$\mathrm{Ind}(L)=0$ but for the sublattice
$M=\{x_2,x_4,x_5,x_6,1_L\}$ of $L$ we have that $\mathrm{Ind}(M)=1$
and thus, $\mathrm{Ind}(M)\nleqslant\mathrm{Ind}(L)$.

\begin{proposition} Let $L$ be a finite lattice and $M$ be a sublattice of $L$ such that
$1_L\in M$, $M\setminus\{1_L\}$ is a down-set and the lattice
${\uparrow}(u^{*}_M\vee u)$ of every element $u$ of $M$ is
isomorphic to the lattice ${\uparrow}(u^{*}_L\vee u)$ of $u$ in $L$.
Then, $\mathrm{Ind}(M)\leqslant\mathrm{Ind}(L)$.
\end{proposition}
{\bf Proof.} If $\mathrm{Ind}(L)=-1$, then, obviously,
$\mathrm{Ind}(M)=-1$. We suppose that $\mathrm{Ind}(L)=k$, where
$k\in \mathbb{N}$, and we will prove that $\mathrm{Ind}(M)\leqslant
k$.

Let $a\in M$ and $v\in M$ with $a\vee v=1_M$, that is $a\vee v=1_L$.
Since $M\subseteq L$, $a,v\in L$. Also, since $\mathrm{Ind}(L)=k$,
there exists $u\in L$ such that $u\leqslant v$, $a\vee u=1_L$ and
$\mathrm{Ind}({\uparrow}(u^{*}_{L} \vee u))\leqslant k-1$. We
consider the following cases:

(a) If $u=1_L$, then clearly $u\in M$ and
$\mathrm{Ind}({\uparrow}(u^{*}_{M}\vee u))=-1=0-1\leqslant k-1$.

(b) Let $u\neq 1_L$. Since $M\setminus\{1_L\}$ is a down set and
$u\leqslant v$, $u\in M$. By assumption since the lattices
${\uparrow}(u^{*}_M\vee u)$ and ${\uparrow}(u^{*}_L\vee u)$ are
isomorphic, we have $\mathrm{Ind}({\uparrow}(u^{*}_{M} \vee
u))=\mathrm{Ind}({\uparrow}(u^{*}_{L} \vee u))\leqslant k-1$.
Therefore, in each case $\mathrm{Ind}(M)\leqslant k$. $\Box$

\section{Comparing the dimension Ind with other lattice notions}\label{sec4}

In this section, we present relations between the large inductive
dimension Ind and the small inductive dimension ind, the Krull
dimension Kdim, the covering dimension dim and the height.

\subsection{The small inductive dimension}\label{subsec2}

\medskip
The notion of the small inductive dimension for finite lattices is
inserted in \cite{GMPS}.

\begin{definition}{\rm \label{ind} Let $L$ be a finite lattice. The \emph{small
inductive dimension}, $\mathrm{ind}$, of $L$ is defined as follows:
\begin{enumerate}
\item $\mathrm{ind}(L)=-1$ if and only if $L=\{0_L\}$.
\item $\mathrm{ind}(L)\leqslant k$, where $k\in \mathbb{N}$, if
for every cover $V$ of $L$, there exists a cover $U$ of $L$ such
that $U$ is a refinement of $V$ and $\mathrm{ind} ({\uparrow}(u^{*}
\vee u))\leqslant k-1$, for every $u\in U$.
\item $\mathrm{ind}(L)=k$, where $k\in \mathbb{N}$, if $\mathrm{ind}(L)\leqslant k$ and $\mathrm{ind}(L)\nleqslant k-1$.
\end{enumerate}}
\end{definition}

The following results will be very useful in our study.

\begin{theorem} \cite{GMPS} \label{minimal cover} Let $L$ be a finite lattice and $k\in\mathbb{N}$. Then $\mathrm{ind}(L)\leqslant k$ if and only if for every minimal cover $V$ of
$L$ we have $\mathrm{ind} ({\uparrow}(v^{*} \vee v))\leqslant k-1$,
for every $v\in V$.
\end{theorem}

However, a modification of Theorem \ref{minimal cover} is given as
follows and can easily be proved following Definition \ref{ind} as
well as Theorem \ref{minimal cover}.

\begin{theorem} \label{minimal cover 2} Let $L$ be a finite lattice and $k\in\mathbb{N}$. Then $\mathrm{ind}(L)\leqslant k$ if and only
if for every $v\in L$ such that $v$ belongs to a minimal cover of
$L$ we have ${\rm ind}({\uparrow}(v^{*} \vee v))\leqslant k-1$.
\end{theorem}

\begin{theorem} \label{mtheo1} If $L$ is a finite lattice, then $\mathrm{ind}(L)\leqslant\mathrm{Ind}(L)$.
\end{theorem}
{\bf Proof.} Clearly, if $\mathrm{Ind}(L)=-1$, the inequality holds.
We suppose that the inequality holds for all finite lattices $M$
with $\mathrm{Ind}(M)\leqslant k-1$. Let $L$ be a finite lattice
such that $\mathrm{Ind}(L)=k$, $w\in L$ and $W=\{w_1,\ldots,w_m\}$
be a minimal cover of $L$ such that $w\in W$. Then $w=w_i$ for some
$i=1,\ldots,m$. By Theorem \ref{minimal cover 2} it suffices to
prove that $\mathrm{ind}({\uparrow}(w_i\vee w_i^{*}))\leqslant k-1.$

Let $a=\displaystyle\bigvee_{j=1,j\neq i}^{m}w_j$ and $v=w_i$ with
$a\vee w_i=1_L$. Since $\mathrm{Ind}(L)=k$, there exists $u\in L$
such that $u\leqslant w_i$, $a\vee u=1_L$ and
$\mathrm{Ind}({\uparrow}(u\vee u^{*}))\leqslant k-1$ and by
inductive assumption, $\mathrm{ind}({\uparrow}(u\vee
u^{*}))\leqslant k-1$. Then the set $U=\{w_1,\ldots,u,\ldots,w_m\}$,
where $U$ is the set $W$ replacing $w_i$ by $u$, is a cover of $L$
which refines $W$. Since $W$ is minimal, $W\subseteq U$. Thus,
$w_i=u$ and then $\mathrm{ind}({\uparrow}(w_i\vee w_i^{*}))\leqslant
k-1$. Thus, $\mathrm{ind}(L)\leqslant k$. $\Box$

\begin{remark}{\rm \label{diff} Generally, the small inductive dimension and the large inductive
dimension are different notions of dimensions for finite lattices.
For example, for the lattice of Figure \ref{fig:4} we have
$\mathrm{Ind}(L)=1$ and $\mathrm{ind}(L)=0$.

Especially, for the large inductive dimension the element which
gives us that $\mathrm{Ind}(L)=1$ is the element $a=x_7$ for which
$v\in\{x_6,x_8,1_L\}$ with $a\vee v=1_L$ and for $u=x_6$ we have
that $\mathrm{Ind}({\uparrow}(x_6\vee x^{*}_6))=0$. For the small
inductive dimension we observe that for all minimal covers $V$ of
$L$ we have $\mathrm{ind}({\uparrow}(v\vee v^{*}))=-1$, for every
$v\in V$. Thus, by Theorem \ref{minimal cover} we have
$\mathrm{ind}(L)=0$.

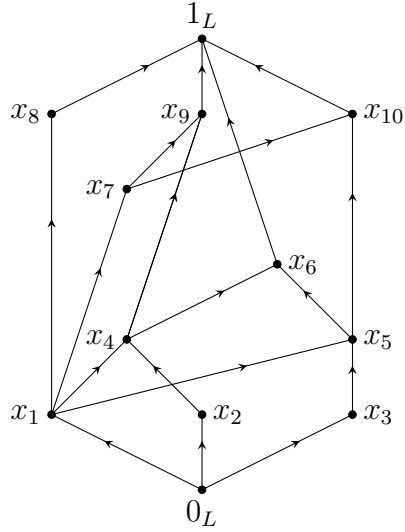
\begin{figure}[H] \tikzstyle arrowstyle=[scale=1]
\tikzstyle
directed=[postaction={decorate,decoration={markings,mark=at position
.65 with {\arrow[arrowstyle]{stealth}}}}]
\begin{center}
\begin{tikzpicture}
\draw[directed] (0,0)--(-2,1); \draw[directed] (0,0)--(0,1);
\draw[directed] (0,0)--(2,1); \draw[directed] (-2,1)--(-1,2);
\draw[directed] (-2,1)--(2,2); \draw[directed] (0,1)--(-1,2);
\draw[directed] (2,1)--(2,2); \draw[directed] (-1,2)--(1,3);
\draw[directed] (2,2)--(1,3); \draw[directed] (-2,1)--(-1,4);
\draw[directed] (-2,1)--(-2,5); \draw[directed] (2,2)--(2,5);
\draw[directed] (-1,2)--(0,5); \draw[directed] (-1,2)--(0,5);
\draw[directed] (-1,4)--(0,5); \draw[directed] (-1,4)--(2,5);
\draw[directed] (1,3)--(0,6); \draw[directed] (-2,5)--(0,6);
\draw[directed] (0,5)--(0,6); \draw[directed] (2,5)--(0,6);
\draw[fill] (0,0) circle [radius=0.05] node [below] {$0_{L}$};
\draw[fill] (-2,1) circle [radius=0.05] node [left] {$x_1$};
\draw[fill] (0,1) circle [radius=0.05] node [right] {$x_2$};
\draw[fill] (2,1) circle [radius=0.05] node [right] {$x_3$};
\draw[fill] (-1,2) circle [radius=0.05] node [left] {$x_4$};
\draw[fill] (2,2) circle [radius=0.05] node [right] {$x_5$};
\draw[fill] (1,3) circle [radius=0.05] node [right] {$x_6$};
\draw[fill] (-1,4) circle [radius=0.05] node [left] {$x_7$};
\draw[fill] (-2,5) circle [radius=0.05] node [left] {$x_8$};
\draw[fill] (0,5) circle [radius=0.05] node [left] {$x_9$};
\draw[fill] (2,5) circle [radius=0.05] node [right] {$x_{10}$};
\draw[fill] (0,6) circle [radius=0.05] node [above] {$1_{L}$};
\end{tikzpicture}
\end{center}
\caption{The finite lattice $(L,\leqslant)$} \label{fig:4}
\end{figure}

We state that the family \[\mathcal{MC}=\{\{x_7,x_8\},
\{x_2,x_3,x_7\}, \{x_2,x_8\}, \{x_3,x_8\}\}\] consists of all the
minimal covers of $L$. The fact that gives us the difference between
the two dimensions $\mathrm{ind}$ and $\mathrm{Ind}$ is the study of
the element $x_6$. For this element we have
$\mathrm{Ind}({\uparrow}(x_6\vee x^{*}_6))=0$ and therefore,
$\mathrm{Ind}(L)=1$. Moreover, we can see that
$\mathrm{ind}({\uparrow}(x_6\vee x^{*}_6))=0$. Thus, if the element
$x_6$ belongs to a minimal cover of $L$, then $\mathrm{ind}(L)=1$
and hence, $\mathrm{ind}(L)=\mathrm{Ind}(L)$. However, we observe
that there is no minimal cover $V$ of $L$ with $x_6\in V$ and thus,
$x_6$ does not make the small inductive dimension equal to $1$. All
minimal covers of $L$ and their elements give us that
$\mathrm{ind}(L)=0$.}
\end{remark}

We can generalize Remark \ref{diff} in order to construct a finite
lattice $L$ with arbitrary large inductive dimension and
$\mathrm{Ind}(L)=\mathrm{ind}(L)+1$.

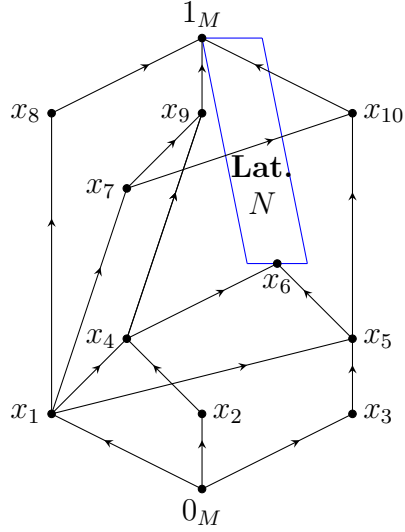
\begin{figure}[H] \tikzstyle arrowstyle=[scale=1]
\tikzstyle
directed=[postaction={decorate,decoration={markings,mark=at position
.65 with {\arrow[arrowstyle]{stealth}}}}]
\begin{center}
\begin{tikzpicture}
\draw[directed] (0,0)--(-2,1); \draw[directed] (0,0)--(0,1);
\draw[directed] (0,0)--(2,1); \draw[directed] (-2,1)--(-1,2);
\draw[directed] (-2,1)--(2,2); \draw[directed] (0,1)--(-1,2);
\draw[directed] (2,1)--(2,2); \draw[directed] (-1,2)--(1,3);
\draw[directed] (2,2)--(1,3); \draw[directed] (-2,1)--(-1,4);
\draw[blue] (1,3)--(0.6,3); \draw[blue] (1,3)--(1.4,3); \draw[blue]
(0.6,3)--(0,6); \draw[blue] (1.4,3)--(0.8,6); \draw[blue]
(0.8,6)--(0,6); \draw[directed] (-2,1)--(-2,5); \draw[directed]
(2,2)--(2,5); \draw[directed] (-1,2)--(0,5); \draw[directed]
(-1,2)--(0,5); \draw[directed] (-1,4)--(0,5); \draw[directed]
(-1,4)--(2,5); \draw[directed] (-2,5)--(0,6); \draw[directed]
(0,5)--(0,6); \draw[directed] (2,5)--(0,6); \draw[fill] (0,0) circle
[radius=0.05] node [below] {$0_{M}$}; \draw[fill] (-2,1) circle
[radius=0.05] node [left] {$x_1$}; \draw[fill] (0,1) circle
[radius=0.05] node [right] {$x_2$}; \draw[fill] (2,1) circle
[radius=0.05] node [right] {$x_3$}; \draw[fill] (-1,2) circle
[radius=0.05] node [left] {$x_4$}; \draw[fill] (2,2) circle
[radius=0.05] node [right] {$x_5$}; \draw[fill] (1,3) circle
[radius=0.05] node [below] {$x_6$}; \draw[fill] (-1,4) circle
[radius=0.05] node [left] {$x_7$}; \draw[fill] (-2,5) circle
[radius=0.05] node [left] {$x_8$}; \draw[fill] (0,5) circle
[radius=0.05] node [left] {$x_9$}; \draw[fill] (2,5) circle
[radius=0.05] node [right] {$x_{10}$}; \draw[fill] (0,6) circle
[radius=0.05] node [above] {$1_{M}$}; \draw (0.8,4.3) node {{\bf
Lat.}}; \draw (0.8,3.8) node {{\bf $N$}};
\end{tikzpicture}
\end{center}
\caption{The finite lattice $(M,\leqslant)$} \label{fig:17}
\end{figure}

\begin{proposition} \label{pr01} Let $k\in\{1,2,\ldots\}$. Then there exists a finite lattice $M$ such that
$\mathrm{ind}(M)=k-1$ and $\mathrm{Ind}(M)=k$.
\end{proposition}
{\bf Proof.} We show the proposition by induction on $k$. Remark
\ref{diff} proves the proposition for the case $k=1$. We consider
the lattice $L$ of Figure \ref{fig:4}. On this lattice we replace
the sublattice $(\{x_6,1_L\},\leqslant)$ with a sublattice $N$ such
that:
\begin{enumerate}
\item
$\mathrm{Ind}(N)=k-1$, where $k\in\{2,3,\ldots\}$,
\item
$x_6=0_N$,
\item
$x_6\leqslant w\leqslant 1_L, \ \mathrm{for} \ \mathrm{every} \ w\in
N$,
\item
the pseudocomplement of every element of $N\setminus\{0_N\}$ is
$0_N$ and
\item
each $u\in L\setminus N\cup\{z\in L:z\leqslant x_6\}$ and each
element $v\in N\setminus\{0_N,1_N\}$ are incomparable (see Figure
\ref{fig:17}).
\end{enumerate}

We state that the existence of the lattice $N$ follows by Theorem
\ref{mth}. We write this new lattice as $M$. Also, for the elements
$a=x_7$ and $v=x_6$, we have $x_7 \vee x_6=1_M$ and
\[\{u \in L : u\leqslant x_6 \ \mathrm{and} \ x_7 \vee u = 1_M\} =
\{x_6\}.\] Since \[\mathrm{Ind}({\uparrow}(x_6 \vee
x_6^*))=\mathrm{Ind}({\uparrow}(x_6\vee
0_M))=\mathrm{Ind}({\uparrow}x_6)=\mathrm{Ind}(N) =k-1,\] we have
$\mathrm{Ind}(M)\geqslant k$. Furthermore, following Definition
\ref{Ind} we can see that $\mathrm{Ind}(M)\leqslant k$.

Moreover, the construction of the lattice $M$ leads us to observe
that the minimal covers of $M$ are the minimal covers of the family
$\mathcal{MC}$ of Remark \ref{diff}, and the minimal covers of $N$.
For each minimal cover $V\in\mathcal{MC}$ we have
$\mathrm{ind}({\uparrow}(v\vee v^{*}))=-1$, for every $v\in V$.
Also, since the construction of $N$ follows the construction given
in \cite[Theorem 4]{GMPS}, we have $\mathrm{ind}(N)=k-1$. Let $U$ be
a minimal cover of $N$. Since $\mathrm{ind}(N)\leqslant k-1$, by
Theorem \ref{minimal cover} we have
\[\mathrm{ind}({\uparrow}(u\vee u^{*}_N))=\mathrm{ind}({\uparrow}(u\vee
x_6))=\mathrm{ind}({\uparrow}u)\leqslant k-2,\] for every $u\in U$,
where $u^{*}_N$ denotes the pseudocomplement of $u$ in the lattice
$N$. Thus, \[\mathrm{ind}({\uparrow}(u\vee
u^{*}_M))=\mathrm{ind}({\uparrow}(u\vee
0_M))=\mathrm{ind}({\uparrow}u)\leqslant k-2,\] for every $u\in U$,
where $u^{*}_M$ denotes the pseudocomplement of $u$ in the lattice
$M$. Therefore, by Theorem \ref{minimal cover}, we have that
$\mathrm{ind}(M)\leqslant k-1$. We prove that $\mathrm{ind}(M)=k-1$.
Since $\mathrm{ind}(N)=k-1$, there exist a minimal cover $W$ of $N$
and $w\in W$ such that \[\mathrm{ind}({\uparrow}(w\vee
w^{*}_N))=\mathrm{ind}({\uparrow}(w\vee
x_6))=\mathrm{ind}({\uparrow}w)=k-2.\] Hence,
$\mathrm{ind}({\uparrow}(w\vee
w^{*}_M))=\mathrm{ind}({\uparrow}(w\vee
0_M))=\mathrm{ind}({\uparrow}w)=k-2$. By Theorem \ref{minimal cover}
we conclude that $\mathrm{ind}(M)=k-1$. $\Box$

\begin{remark}{\rm \label{gap} We observe that the ``gap" between the
dimensions $\mathrm{ind}$ and $\mathrm{Ind}$ can be arbitrarily
large. For that, we consider the lattice $(Q,\leqslant)$ of Figure
\ref{fig:18}. For this lattice we have that $\mathrm{Ind}(Q)=3$ and
$\mathrm{ind}(Q)=0$.

Firstly, we describe the construction of the lattice
$(Q,\leqslant)$. We consider the lattice $M$ of Figure \ref{fig:17},
where $N=\{1_M,x_6\}\cup\{w_1,\ldots,w_6\}$ is the sublattice given
in Theorem \ref{mth} for $k=2$, that is $\mathrm{Ind}(N)=2$. On this
lattice $M$, we add the following elements:

\noindent (I) for each $w\in N\setminus\{1_M,x_6\}$ we assign an
element $x_{w}$ such that:

(1) $0_M<x_w<w$,

(2) for every $w_1,w_2 \in N\setminus\{1_M,x_6\}$ with $w_1\neq
w_2$, $x_{w_{1}}\neq x_{w_{2}}$,

(3) $x_w$ is incomparable with any other element of $M$,

\noindent (II) an element $s$ such that:

(1) $0_M<s<x_w$, for every $x_w$,

\noindent (III) elements $t,r,u$ as in Figure \ref{fig:18} in order
to succeed that the constructed poset is a lattice.

\begin{figure}[H] \tikzstyle arrowstyle=[scale=1]
\tikzstyle
directed=[postaction={decorate,decoration={markings,mark=at position
.65 with {\arrow[arrowstyle]{stealth}}}}]
\begin{center}
\begin{tikzpicture}
\draw[directed] (3,0)--(0,1); \draw[directed] (3,0)--(7,1);
\draw[directed] (3,0)--(-2,1); \draw[directed] (-2,1)--(-1,2);
\draw[directed] (0,1)--(-1,2); \draw[directed] (7,1)--(7,2);
\draw[directed] (-2,1)--(7,2); \draw[directed] (4,4)--(4,5);
\draw[directed] (4,5)--(3,6); \draw[directed] (4,5)--(6,8);
\draw[directed] (3,6)--(3,7); \draw[directed] (3,7)--(2,8);
\draw[directed] (2,8)--(3,9); \draw[directed] (3,7)--(4,8);
\draw[directed] (4,8)--(3,9); \draw[directed] (6,8)--(3,9);
\draw[directed] (3,4)--(4,5); \draw[directed] (2,4)--(3,6);
\draw[directed] (1,4)--(3,7); \draw[directed] (0,4)--(2,8);
\draw[directed] (5,4)--(4,8); \draw[directed] (6,4)--(6,8);
\draw[directed] (-2,1)--(-2,8); \draw[directed] (-2,1)--(-1.5,6.4);
\draw[directed] (-1,2)--(0,8); \draw[directed] (-1.5,6.4)--(0,8);
\draw[directed] (-1.5,6.4)--(7,8); \draw[directed] (-2,8)--(3,9);
\draw[directed] (0,8)--(3,9); \draw[directed] (7,8)--(3,9);
\draw[directed] (3,0)--(4,1); \draw[directed] (4,1)--(4,2);
\draw[directed] (4,1)--(0,4);  \draw[directed] (4,1)--(1,4);
\draw[directed] (4,1)--(2,4);  \draw[directed] (4,1)--(3,4);
\draw[directed] (4,1)--(5,4);  \draw[directed] (4,1)--(6,4);
\draw[directed] (-2,1)--(4,2);  \draw[directed] (4,2)--(4,4);
\draw[directed] (4,2)--(-1.5,6.4); \draw[directed] (4,2)--(4,3);
\draw[directed] (-1,2)--(4,3); \draw[directed] (4,3)--(0,8);
\draw[directed] (4,2)--(7,3); \draw[directed] (7,2)--(7,3);
\draw[directed] (7,3)--(4,4); \draw[directed] (7,3)--(7,8);
\draw[fill] (3,0) circle [radius=0.05] node [below] {$0_{Q}$};
\draw[fill] (0,1) circle [radius=0.05] node [right] {$x_2$};
\draw[fill] (7,1) circle [radius=0.05] node [right] {$x_3$};
\draw[fill] (7,2) circle [radius=0.05] node [right] {$x_5$};
\draw[fill] (-1,2) circle [radius=0.05] node [left] {$x_4$};
\draw[fill] (-2,1) circle [radius=0.05] node [left] {$x_1$};
\draw[fill] (4,4) circle [radius=0.05] node [right] {$x_6$};
\draw[fill] (4,5) circle [radius=0.05] node [right] {$w_1$};
\draw[fill] (3,6) circle [radius=0.05] node [right] {$w_2$};
\draw[fill] (3,7) circle [radius=0.05] node [left] {$w_3$};
\draw[fill] (2,8) circle [radius=0.05] node [left] {$w_4$};
\draw[fill] (4,8) circle [radius=0.05] node [right] {$w_5$};
\draw[fill] (6,8) circle [radius=0.05] node [left] {$w_6$};
\draw[fill] (0,4) circle [radius=0.05] node [left] {$x_{w_4}$};
\draw[fill] (1,4) circle [radius=0.05] node [left] {$x_{w_3}$};
\draw[fill] (2,4) circle [radius=0.05] node [left] {$x_{w_2}$};
\draw[fill] (3,4) circle [radius=0.05] node [left] {$x_{w_1}$};
\draw[fill] (5,4) circle [radius=0.05] node [right] {$x_{w_5}$};
\draw[fill] (6,4) circle [radius=0.05] node [right] {$x_{w_6}$};
\draw[fill] (-2,8) circle [radius=0.05] node [left] {$x_8$};
\draw[fill] (-1.5,6.4) circle [radius=0.05] node [left] {$x_7$};
\draw[fill] (0,8) circle [radius=0.05] node [left] {$x_9$};
\draw[fill] (7,8) circle [radius=0.05] node [right] {$x_{10}$};
\draw[fill] (4,1) circle [radius=0.05] node [right] {$s$};
\draw[fill] (4,2) circle [radius=0.05] node [right] {$t$};
\draw[fill] (4,3) circle [radius=0.05] node [right] {$r$};
\draw[fill] (7,3) circle [radius=0.05] node [right] {$u$};
\draw[fill] (3,9) circle [radius=0.05] node [above] {$1_Q$};
\end{tikzpicture}
\end{center}
\caption{The finite lattice $(Q,\leqslant)$} \label{fig:18}
\end{figure}
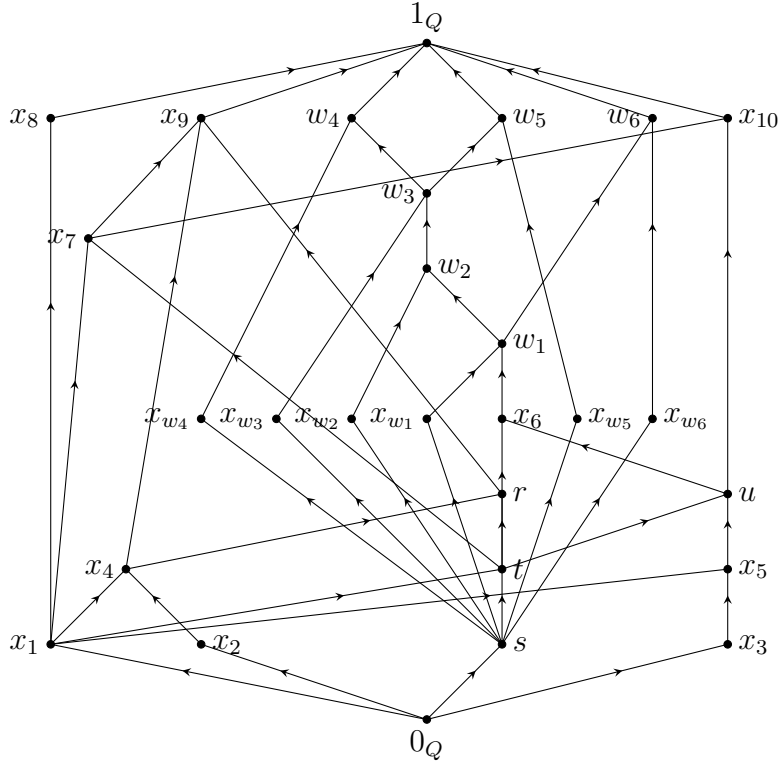

We write this new lattice as $Q$. Using the same argument as in
Proposition \ref{pr01} we have $\mathrm{Ind}(Q)=3$, considering the
elements $a=x_7$ and $v=x_6$. Also, the minimal covers of $Q$ are
the family
\[\mathcal{MC}^{\prime}=\{\{s,x_8\}, \{x_2,x_3,x_7\}, \{x_2,x_8\}, \{x_3,x_8\}\}\] and the minimal covers
consisting of elements of the set $\{x_{w_{1}},\ldots,x_{w_{6}}\}$.
For each minimal cover $V\in\mathcal{MC}^{\prime}$ we have
$\mathrm{ind}({\uparrow}(v\vee v^{*}))=-1$, for every $v\in V$.
Also, for each $x_w\in Q$ we have
\[\mathrm{ind}({\uparrow}(x_w\vee x_{w}^{*}))=\mathrm{ind}({\uparrow}(x_w\vee 1_Q))=\mathrm{ind}({\uparrow}1_Q)=\mathrm{ind}(\{1_Q\})=-1.\] Therefore,
$\mathrm{ind}(Q)=0$.}
\end{remark}

Now, we can generalize Remark \ref{gap} in order to construct a
finite lattice $P$ for which $\mathrm{Ind}(P)-\mathrm{ind}(P)$ is
arbitrarily large.

\begin{figure}[H] \tikzstyle arrowstyle=[scale=1]
\tikzstyle
directed=[postaction={decorate,decoration={markings,mark=at position
.65 with {\arrow[arrowstyle]{stealth}}}}]
\begin{center}
\begin{tikzpicture}
\draw[directed] (0,0)--(-2,1); \draw[directed] (0,0)--(0,1);
\draw[directed] (0,0)--(2,1); \draw[directed] (-2,1)--(-1,2);
\draw[directed] (-2,1)--(2,2); \draw[directed] (0,1)--(-1,2);
\draw[directed] (2,1)--(2,2); \draw[directed] (-2,1)--(-1,4);
\draw[blue] (1,3)--(0.6,3); \draw[blue] (1,3)--(1.4,3); \draw[blue]
(0.6,3)--(0,6); \draw[blue] (1.4,3)--(0.8,6); \draw[blue]
(0.8,6)--(0,6); \draw[directed] (-2,1)--(-2,5); \draw[directed]
(-1,2)--(0,5); \draw[directed] (-1,2)--(0,5); \draw[directed]
(-1,4)--(0,5); \draw[directed] (-1,4)--(2,5); \draw[directed]
(-2,5)--(0,6); \draw[directed] (0,5)--(0,6); \draw[directed]
(2,5)--(0,6); \draw[directed] (0,0)--(-4,1); \draw[directed]
(-4,1)--(-4,2); \draw[directed] (-2,1)--(-4,2); \draw[directed]
(-4,2)--(-0.3,2.5); \draw[directed] (-1,2)--(-0.3,2.5);
\draw[directed] (-0.3,2.5)--(1,3); \draw[directed] (-4,2)--(2,2.5);
\draw[directed] (2,2)--(2,2.5); \draw[directed] (2,2.5)--(1,3);
\draw[directed] (2,2.5)--(2,5); \draw[directed] (-0.3,2.5)--(0,5);
\draw[directed] (-4,2)--(-1,4); \draw[blue] (-8.5,3)--(-4,3);
\draw[blue] (-4,3)--(-4,4); \draw[blue] (-4,4)--(-8.5,4);
\draw[blue] (-8.5,4)--(-8.5,3); \draw[blue] (-4,1)--(-6,3);
\draw[blue] (-4,3)--(1,4); \draw[blue] (-4,4)--(0.8,5); \draw[fill]
(0,0) circle [radius=0.05] node [below] {$0_{P}$}; \draw[fill]
(-2,1) circle [radius=0.05] node [left] {$x_1$}; \draw[fill] (0,1)
circle [radius=0.05] node [right] {$x_2$}; \draw[fill] (2,1) circle
[radius=0.05] node [right] {$x_3$}; \draw[fill] (-1,2) circle
[radius=0.05] node [left] {$x_4$}; \draw[fill] (2,2) circle
[radius=0.05] node [right] {$x_5$}; \draw[fill] (1,3) circle
[radius=0.05] node [below] {$x_6$}; \draw[fill] (-1,4) circle
[radius=0.05] node [left] {$x_7$}; \draw[fill] (-2,5) circle
[radius=0.05] node [left] {$x_8$}; \draw[fill] (0,5) circle
[radius=0.05] node [left] {$x_9$}; \draw[fill] (2,5) circle
[radius=0.05] node [right] {$x_{10}$}; \draw[fill] (0,6) circle
[radius=0.05] node [above] {$1_{P}$}; \draw[fill] (-4,1) circle
[radius=0.05] node [left] {$s$}; \draw[fill] (-4,2) circle
[radius=0.05] node [left] {$t$}; \draw[fill] (-0.3,2.5) circle
[radius=0.05] node [right] {$r$}; \draw[fill] (2,2.5) circle
[radius=0.05] node [right] {$u$}; \draw (0.8,4.3) node {{\bf Lat.}};
\draw (0.8,3.8) node {{\bf $N$}}; \draw (-6.25,3.5) node
{$\{x_w:w\in N\setminus\{1_P,x_6\}\}$};
\end{tikzpicture}
\end{center}
\caption{The finite lattice $(P,\leqslant)$} \label{fig:19}
\end{figure}
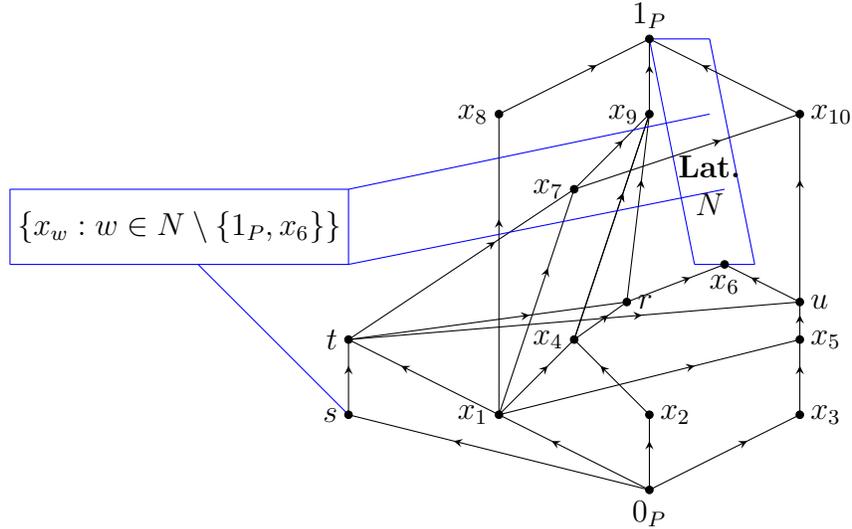

\begin{proposition} Let $k\in\{1,2,\ldots\}$. Then there exists a finite lattice $P$ such that
$\mathrm{ind}(P)=0$ and $\mathrm{Ind}(P)=k$.
\end{proposition}
{\bf Proof.} We shall prove the proposition by induction on $k$. If
$k=1$, then it follows by Proposition \ref{pr01} that we get the
desired result. Suppose $k\geqslant 2$ and we have a finite lattice
$N$ such that $\mathrm{ind}(N)=0$ and $\mathrm{Ind}(N)=k-1$. We
consider the lattice $M$ of Figure \ref{fig:17}. On this lattice
$M$, we add the following elements (Figure \ref{fig:19}):

\noindent (I) for each $w\in N\setminus\{1_M,x_6\}$ we assign an
element $x_{w}$ such that:

(1) $0_M<x_w<w$,

(2) for every $w_1,w_2 \in N\setminus\{1_M,x_6\}$ with $w_1\neq
w_2$, $x_{w_{1}}\neq x_{w_{2}}$,

(3) $x_w$ is incomparable with any other element of $M$,

\noindent (II) an element $s$ such that:

(1) $0_M<s<x_w$, for every $x_w$,

\noindent (III) elements $t,r,u$ as in Figure \ref{fig:19} in order
to succeed that the constructed poset is a lattice.

We write this new lattice as $P$. Using the same argument as in
Proposition \ref{pr01} we have $\mathrm{Ind}(P)=k$. Also, the
minimal covers of $P$ are the family $\mathcal{MC}^{\prime}$ of
Remark \ref{gap} and the minimal covers consisting of elements of
the set $\{x_w:w\in N\setminus\{1_P,x_6\}\}$. For each minimal cover
$V\in\mathcal{MC}^{\prime}$ we have $\mathrm{ind}({\uparrow}(v\vee
v^{*}))=-1$, for every $v\in V$. Also, for each $x_w\in P$ we have
\[\mathrm{ind}({\uparrow}(x_w\vee x_{w}^{*}))=\mathrm{ind}({\uparrow}(x_w\vee 1_P))=\mathrm{ind}({\uparrow}1_P)=\mathrm{ind}(\{1_P\})=-1.\] Therefore,
$\mathrm{ind}(P)=0$. $\Box$

\medskip
However, we observe that if we add some conditions for finite
lattices, then we can have the equality of these two dimensions.

\begin{lemma}\label{l} Let $L$ be a finite lattice and $M={\uparrow}x$,
where $x\in L$. If for every minimal cover $V$ of $L$ we have
$|V|\leqslant 2$, then for every minimal cover $W$ of $M$ we have
also $|W|\leqslant 2$.
\end{lemma}
{\bf Proof.} Let $W$ be a minimal cover of $M$. We shall prove that
$|W|\leqslant 2$. We have that $W$ is a cover of $L$.

\noindent (1) If $W$ is a minimal cover of $L$, then by assumption
$|W|\leqslant 2$.

\noindent (2) We suppose that $W$ is not a minimal cover of $L$.
Then we consider the following cases:

(a) If $1_L\in W$ that is $1_M\in W$ (as $1_L=1_M$), then
$W=\{1_M\}$ is the only minimal cover of $M$ for which $|W|=1$.

(b) We state that $1_L\notin W$. By \cite[Lemma 3.7]{DGMM}, there
exists a minimal cover $R$ of $L$ which refines $W$ and $1_L\notin
R$. By assumption we have that $|R|=\{c,d\}$. Then $c\leqslant a$
and $d\leqslant b$, for some $a,b\in W$. Moreover, $a\vee
b=1_L=1_M$. Indeed, if $a\vee b=r$, where $r\neq 1_L$, then the set
$\{c,d\}$ is not a cover of $L$, which is a contradiction. Thus,
$\{a,b\}$ is a cover of $M$, which refines $W$. Since $W$ is a
minimal cover of $M$ and $a,b\in W$, we have that $W=\{a,b\}$, that
is $|W|=2$.

Therefore, in each case we have that $|W|\leqslant 2$. $\Box$

\begin{proposition}Let $L$ be a finite lattice. If for every minimal cover $V$ of $L$ we have $|V|\leqslant
2$, then $\mathrm{ind}(L)=\mathrm{Ind}(L)$.
\end{proposition}
{\bf Proof.} By Theorem \ref{mtheo1} we have that
$\mathrm{ind}(L)\leqslant\mathrm{Ind}(L)$. Thus, it suffices to
prove that $\mathrm{Ind}(L)\leqslant\mathrm{ind}(L)$. Clearly, if
$\mathrm{ind}(L)=-1$, the inequality holds. We suppose that the
proposition holds for all finite lattices $M$ with
$\mathrm{ind}(M)\leqslant k-1$. Let $L$ be a finite lattice such
that $\mathrm{ind}(L)=k$. We shall prove that
$\mathrm{Ind}(L)\leqslant k$.

(1) For $a=1_L$ and any $v\in L$ with $a\vee v=1_L$, there exists
the element $u=0_L$ such that $u\leqslant v$, $a\vee u=1_L$ and
$\mathrm{Ind}({\uparrow}(u\vee u^{*}))=-1\leqslant k-1$.

(2) For $a=0_L$ and $v=1_L$ with $a\vee v=1_L$, there exists the
element $u=1_L$ such that $u\leqslant v$, $a\vee u=1_L$ and
$\mathrm{Ind}({\uparrow}(u\vee u^{*}))=-1\leqslant k-1$.

(3) Let $a\in L\setminus\{0_L,1_L\}$ and $v\in L$ with $a\vee
v=1_L$. If $v=1_L$, we have the trivial fact
$\mathrm{Ind}({\uparrow}(v\vee v^{*}))=-1$. Therefore, we suppose
that $v\neq 1_L$ and we study the following cases.

(a) If $\{a,v\}$ is a minimal cover of $L$, then since
$\mathrm{ind}(L)=k$, by Theorem \ref{minimal cover} we have
$\mathrm{ind}({\uparrow}(v\vee v^{*}))\leqslant k-1$. Moreover, by
Lemma \ref{l}, the assumption of the proposition also holds for the
lattice ${\uparrow}(v\vee v^{*})$. Therefore, by induction, we have
that $\mathrm{Ind}({\uparrow}(v\vee v^{*}))\leqslant k-1$. Thus, in
Definition \ref{Ind}, if $u=v$, then $u$ is an element of $L$ for
which $a\vee u=1_L$ and $\mathrm{Ind}({\uparrow}(u\vee
u^{*}))\leqslant k-1$.

(b) We suppose that $\{a,v\}$ is not a minimal cover of $L$. We
state that $1_L\notin \{a,v\}$. By \cite[Lemma 3.7]{DGMM}, there
exists a minimal cover $R$ of $L$ which refines $\{a,v\}$ and
$1_L\notin R$. By assumption we have that $R=\{c,d\}$. Then
$c\leqslant a$ and $d\leqslant v$. Also, since $\mathrm{ind}(L)=k$,
by Theorem \ref{minimal cover} we have
$\mathrm{ind}({\uparrow}(d\vee d^{*}))\leqslant k-1$ and by
induction (as Lemma \ref{l} shows that the assumption of the
proposition also holds for the lattice ${\uparrow}(d\vee d^{*})$),
we have that \[\mathrm{Ind}({\uparrow}(d\vee d^{*}))\leqslant k-1.\]
Moreover, we observe that $a\vee d=1_L$. Indeed, if $a\vee d=r$,
where $r\neq 1_L$, then the set $\{c,d\}$ is not a cover of $L$,
which is a contradiction. Thus, in Definition \ref{Ind}, if $u=d$,
then $u$ is an element of $L$ for which $u\leqslant v$, $a\vee
u=1_L$ and $\mathrm{Ind}({\uparrow}(u\vee u^{*}))\leqslant k-1$.

Therefore, in each case we have that $\mathrm{Ind}(L)\leqslant k$.
$\Box$

\subsection{The Krull dimension}\label{subsec3}

\medskip
A non-empty subset $F$ of a lattice $(L,\leqslant)$ is called
\emph{filter} if $F$ has the following properties:\\
(1) $F\neq L$,\\
(2) If $x\in F$ and $x\leqslant y$, then $y\in F$,\\
(3) If $x,y\in F$, then $x\wedge y\in F$.

A filter $F$ is called \emph{prime} if for every $x,y\in L$ with
$x\vee y\in F$, we have $x\in F$ or $y\in F$. The set of all prime
filters of a lattice $L$ is usually denoted by $\mathcal{PF}(L)$.

\begin{definition}{\rm \cite{V} If $\mathcal{PF}(L)\neq\emptyset$, then the
\emph{Krull dimension} of $L$ is defined as follows:
\[\mathrm{Kdim}(L)=\sup\{k:\mathrm{there \ exist \ prime \ filters} \ F_{0}\subset F_{1}\subset\cdots\subset
F_{k}\}.\]}
\end{definition}

\noindent Also, in the study \cite{DGMS}, the Krull dimension is
studied through join prime elements and matrices. An element $x\in
L$ is said to be \emph{join prime} if it is non-zero and the
inequality $x\leqslant a\vee b$ implies $x\leqslant a$ or
$x\leqslant b$, for all $a,b\in L$.

\begin{example}{\rm We show by examples that the dimensions $\mathrm{Ind}$ and $\mathrm{Kdim}$ for
finite lattices are in general different.

\medskip
\noindent (1) The large inductive dimension of the pentagon (see
Figure \ref{fig:5}) is equal to 0, but since it is not a Boolean
Algebra, its Krull dimension is not equal to 0 (see \cite{V}).

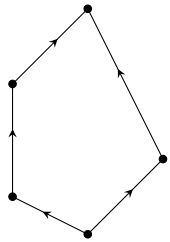
\begin{figure}[H]
\tikzstyle arrowstyle=[scale=1] \tikzstyle
directed=[postaction={decorate,decoration={markings,mark=at position
.6 with {\arrow[arrowstyle]{stealth}}}}]
\begin{center}
\begin{tikzpicture}
\draw[directed] (0,0)--(-1,0.5); \draw[directed] (0,0)--(1,1);
\draw[directed] (-1,0.5)--(-1,2); \draw[directed] (-1,2)--(0,3);
\draw[directed] (1,1)--(0,3); \draw[fill] (0,0) circle
[radius=0.05]; \draw[fill] (-1,0.5) circle [radius=0.05];
\draw[fill] (1,1) circle [radius=0.05]; \draw[fill] (-1,2) circle
[radius=0.05]; \draw[fill] (0,3) circle [radius=0.05];
\end{tikzpicture}
\end{center}
\caption{The pentagon} \label{fig:5}
\end{figure}

\noindent (2) We consider the finite lattice $(L,\leqslant)$,
represented by the diagram of Figure \ref{fig:6}.

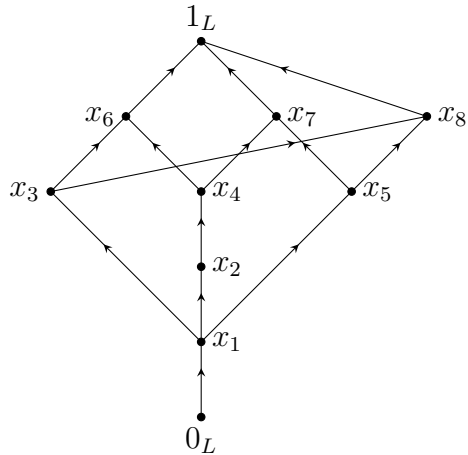
\begin{figure}[H] \tikzstyle arrowstyle=[scale=1]
\tikzstyle
directed=[postaction={decorate,decoration={markings,mark=at position
.65 with {\arrow[arrowstyle]{stealth}}}}]
\begin{center}
\begin{tikzpicture}
\draw[directed] (0,0)--(0,1); \draw[directed] (0,1)--(0,2);
\draw[directed] (0,1)--(-2,3); \draw[directed] (0,1)--(2,3);
\draw[directed] (0,2)--(0,3); \draw[directed] (-2,3)--(-1,4);
\draw[directed] (-2,3)--(3,4); \draw[directed] (0,3)--(-1,4);
\draw[directed] (0,3)--(1,4); \draw[directed] (2,3)--(1,4);
\draw[directed] (2,3)--(3,4); \draw[directed] (-1,4)--(0,5);
\draw[directed] (1,4)--(0,5); \draw[directed] (3,4)--(0,5);
\draw[fill] (0,0) circle [radius=0.05] node [below] {$0_{L}$};
\draw[fill] (0,1) circle [radius=0.05] node [right] {$x_1$};
\draw[fill] (0,2) circle [radius=0.05] node [right] {$x_2$};
\draw[fill] (-2,3) circle [radius=0.05] node [left] {$x_3$};
\draw[fill] (0,3) circle [radius=0.05] node [right] {$x_4$};
\draw[fill] (2,3) circle [radius=0.05] node [right] {$x_5$};
\draw[fill] (-1,4) circle [radius=0.05] node [left] {$x_6$};
\draw[fill] (1,4) circle [radius=0.05] node [right] {$x_7$};
\draw[fill] (3,4) circle [radius=0.05] node [right] {$x_8$};
\draw[fill] (0,5) circle [radius=0.05] node [above] {$1_{L}$};
\end{tikzpicture}
\end{center}
\caption{The finite lattice $(L,\leqslant)$} \label{fig:6}
\end{figure}

We have that $\mathrm{Ind}(L)=2$ and $\mathrm{Kdim}(L)=1$. Indeed,
the elements which give us that $\mathrm{Ind}(L)=2$ are $a=x_{8}$
and $v=x_2$ for which $a \vee v=1_L$. Then $$\{u \in L : u \leqslant
x_2 \ \mathrm{and} \ a \vee u = 1_L\} = \{x_2\}$$ and
$\mathrm{Ind}({\uparrow}(x_2 \vee
x_2^*))=\mathrm{Ind}({\uparrow}x_2)=1$. Hence, $\mathrm{Ind}
(L)\geqslant 2$. It is easy to show that $\mathrm{Ind}(L)\leqslant
2$. Moreover, the join prime elements of $L$ are the elements $x_1$,
$x_2$, $x_3$ and $x_5$ and thus, by Proposition 4.5 (5) of
\cite{DGMS} we have that $\mathrm{Kdim}(L)=1$.}
\end{example}

\subsection{The covering dimension}\label{subsec4}

\medskip
Let $L$ be a finite lattice. The \emph{order} of a subset $C$ of
$L$, denoted by $\mathrm{ord}(C)$, is defined to be $k$, where $k\in
\mathbb{N}$ if and only if the infimum of any $k+2$ distinct
elements of $C$ is $0_L$ and there exist $k+1$ distinct elements of
$C$ whose infimum is not $0_L$.

\begin{definition}{\rm \cite{DGMM} The function $\mathrm{dim}$, called \emph{covering dimension}, with
domain the class of all finite lattices and range the set
$\mathbb{N}$ is defined as follows:
\begin{enumerate}
\item $\mathrm{dim}(L)\leqslant k$, where $k\in\mathbb{N}$, if
and only if for every cover $C$ of $L$, there exists a cover $R$ of
$L$, refinement of $C$ with $\mathrm{ord}(R)\leqslant k$.
\item $\mathrm{dim}(L)=k$, where $k\in \mathbb{N}$, if $\mathrm{dim(L)}\leqslant k$ and $\mathrm{dim}(L)\nleqslant k-1$.
\end{enumerate}}
\end{definition}

In Theorem 3.8 of \cite{DGMM}, it was proved that
$\mathrm{dim}(L)\leqslant k$, where $k\in\mathbb{N}$, if and only if
for every minimal cover $V$ of $L$ we have $\mathrm{ord}(V)\leqslant
k$. That is, $${\rm dim}(L)=\max\{{\rm ord}(V):V \ \mbox{is minimal
cover of} \ L\}$$ (see Corollary 3.9 of \cite{DGMM}).

\begin{remark}{\rm We show by examples that the dimensions $\mathrm{Ind}$ and $\mathrm{dim}$ for
finite lattices are in general different.

\medskip
\noindent (1) We consider the lattice of Figure \ref{fig:7}. We have
that $\mathrm{Ind}(L)=3$ and $\mathrm{dim}(L)=2$.

Indeed, the elements which give us that $\mathrm{Ind}(L)=3$ are
$a=x_{11}$ and $v = x_2$ for which $a \vee v = 1_L$. Then $$\{u \in
L : u\leqslant x_2 \ \mathrm{and} \ a \vee u = 1_L\} = \{x_2\}$$ and
$\mathrm{Ind}({\uparrow}(x_2 \vee
x_2^*))=\mathrm{Ind}({\uparrow}x_2) = 2$. Hence,
$\mathrm{Ind}(L)\geqslant 3$. It is easy to show that
$\mathrm{Ind}(L)\leqslant 3$. Moreover, for all minimal covers $V$
of $L$, that is $V\in\{\{x_2,x_{11}\}, \{x_4,x_5,x_7\}\}$, we have
$\mathrm{ord}(V)\leqslant 2$. Thus, $\mathrm{dim}(L)=2$.

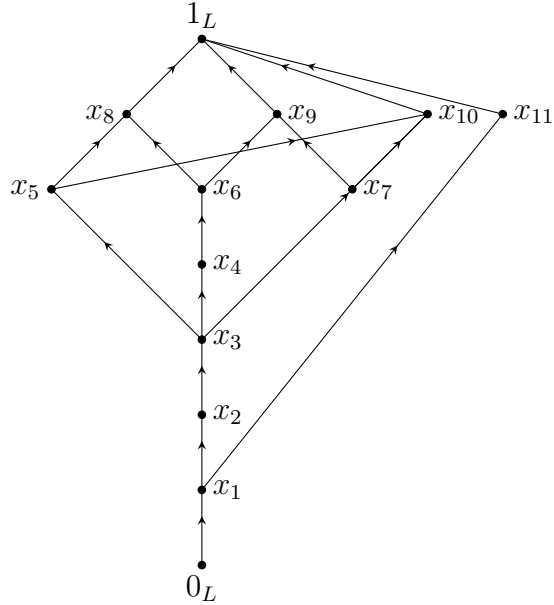
\begin{figure}[H] \tikzstyle arrowstyle=[scale=1]
\tikzstyle
directed=[postaction={decorate,decoration={markings,mark=at position
.65 with {\arrow[arrowstyle]{stealth}}}}]
\begin{center}
\begin{tikzpicture}
\draw[directed] (0,0)--(0,1); \draw[directed] (0,1)--(0,2);
\draw[directed] (0,2)--(0,3); \draw[directed] (0,3)--(0,4);
\draw[directed] (0,4)--(0,5); \draw[directed] (0,3)--(-2,5);
\draw[directed] (0,3)--(3,6); \draw[directed] (-2,5)--(-1,6);
\draw[directed] (-2,5)--(3,6); \draw[directed] (0,5)--(-1,6);
\draw[directed] (0,5)--(1,6); \draw[directed] (2,5)--(1,6);
\draw[directed] (2,5)--(3,6); \draw[directed] (0,1)--(4,6);
\draw[directed] (-1,6)--(0,7); \draw[directed] (1,6)--(0,7);
\draw[directed] (3,6)--(0,7); \draw[directed] (4,6)--(0,7);
\draw[fill] (0,0) circle [radius=0.05] node [below] {$0_{L}$};
\draw[fill] (0,1) circle [radius=0.05] node [right] {$x_1$};
\draw[fill] (0,2) circle [radius=0.05] node [right] {$x_2$};
\draw[fill] (0,3) circle [radius=0.05] node [right] {$x_3$};
\draw[fill] (0,4) circle [radius=0.05] node [right] {$x_4$};
\draw[fill] (-2,5) circle [radius=0.05] node [left] {$x_5$};
\draw[fill] (0,5) circle [radius=0.05] node [right] {$x_6$};
\draw[fill] (2,5) circle [radius=0.05] node [right] {$x_7$};
\draw[fill] (-1,6) circle [radius=0.05] node [left] {$x_8$};
\draw[fill] (1,6) circle [radius=0.05] node [right] {$x_9$};
\draw[fill] (3,6) circle [radius=0.05] node [right] {$x_{10}$};
\draw[fill] (4,6) circle [radius=0.05] node [right] {$x_{11}$};
\draw[fill] (0,7) circle [radius=0.05] node [above] {$1_{L}$};
\end{tikzpicture}
\end{center}
\caption{The finite lattice $(L,\leqslant)$} \label{fig:7}
\end{figure}

\noindent (2) For the lattice $L$ of Figure \ref{fig:8} we have
$\mathrm{dim}(L)=2$ and $\mathrm{Ind}(L)=1$. Indeed, let $a \in
L\setminus\{0_L, 1_L\}$ and $v \in L$ with $a \vee v = 1_L$. It is
clear that $a \vee x_1=a \ne 1_L$. Thus, for each $u \in L$ with
$u\leqslant v$ and $a \vee u = 1_L$, we have $u \ne x_1$. Further,
for each $i = 2, \dots , 7$ $x_i^* = 0_L$ and hence ${\uparrow}(x_i
\vee x_i^*)={\uparrow}x_i$. It is easy to show that
$\mathrm{Ind}({\uparrow}x_i)=0$. Thus, we have
$\mathrm{Ind}(L)\geqslant 1$. Further, we notice that $x_2 \vee x_5
= 1_L$, $$\{u \in L : u\leqslant x_5 \ \mathrm{and} \ x_2 \vee u =
1_L\} = \{x_5\}$$ and $\mathrm{Ind}({\uparrow}x_5)=0$. Hence
$\mathrm{Ind}(L)=1$. For the covering dimension, the only minimal
cover of $L$ is the set $V=\{x_2,x_3,x_4\}$ for which
$\mathrm{ord}(V)=2$. Thus, $\mathrm{dim}(L)=2$.

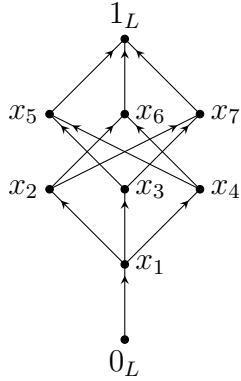
\begin{figure}[H]
\tikzstyle arrowstyle=[scale=1] \tikzstyle
directed=[postaction={decorate,decoration={markings,mark=at position
.85 with {\arrow[arrowstyle]{stealth}}}}]
\begin{center}
\begin{tikzpicture}
\draw[directed] (0,0)--(0,1); \draw[directed] (0,1)--(-1,2);
\draw[directed] (0,1)--(0,2); \draw[directed] (0,1)--(1,2);
\draw[directed] (-1,2)--(0,3); \draw[directed] (-1,2)--(1,3);
\draw[directed] (0,2)--(-1,3); \draw[directed] (0,2)--(1,3);
\draw[directed] (-1,3)--(0,4); \draw[directed] (0,3)--(0,4);
\draw[directed] (1,3)--(0,4); \draw[directed] (1,2)--(-1,3);
\draw[directed] (1,2)--(0,3); \draw[fill] (0,0) circle [radius=0.05]
node [below] {$0_{L}$}; \draw[fill] (0,1) circle [radius=0.05] node
[right] {$x_1$}; \draw[fill] (-1,2) circle [radius=0.05] node [left]
{$x_2$}; \draw[fill] (0,2) circle [radius=0.05] node [right]
{$x_3$}; \draw[fill] (1,2) circle [radius=0.05] node [right]
{$x_4$}; \draw[fill] (-1,3) circle [radius=0.05] node [left]
{$x_5$}; \draw[fill] (0,3) circle [radius=0.05] node [right]
{$x_6$}; \draw[fill] (1,3) circle [radius=0.05] node [right]
{$x_7$}; \draw[fill] (0,4) circle [radius=0.05] node [above]
{$1_L$};
\end{tikzpicture}
\end{center}
\caption{The finite lattice $(L,\leqslant)$} \label{fig:8}
\end{figure}
}
\end{remark}

\subsection{The height}\label{subsec5}

\medskip
\begin{definition}{\rm The {\it height} of a finite lattice $L$, denoted by $h(L)$, is defined as follows:
$$h(L)=\max\{k: \mbox{there exist} \ a_0,\ldots,a_k\in L \ \mbox{such that} \
0_L=a_0<\ldots<a_k=1_L\}.$$}
\end{definition}

\begin{example} (1) For the lattice $L_1$ of Figure \ref{fig:1} we have that ${\rm Ind}(L_1)=0$ and
$h(L_1)=2$.

\medskip
\noindent (2) For the totally ordered set
$L=\{x_1,\ldots,x_n,0_L,1_L\}$ we have that $h(L)=n+1$ and ${\rm
Ind}(L)=0$.
\end{example}

\begin{proposition} \label{height} For every finite lattice $L$, ${\rm Ind}(L)<h(L)$.
\end{proposition}
{\bf Proof.} We prove the inequality by induction on $h(L)$. Clearly
if $h(L)=0$, then $L=\{1_L\}$ and thus, ${\rm Ind}(L)=-1$, proving
that the relation of the proposition holds. We suppose that the
inequality holds for all finite lattices $M$ with $h(M)<k$, for
$k>1$, and we shall prove it for $k$. Let $L$ be a finite lattice
with $h(L)=k$, $v\in L$ such that $v$ belongs to a minimal cover $V$
of $L$. By Theorem \ref{minimal cover 2} it suffices to prove that
${\rm Ind}({\uparrow}(v^{*}\vee v))\leqslant k-2$. Since
${\uparrow}(v^{*}\vee v)\neq L$, we have that
$h({\uparrow}(v^{*}\vee v))<h(L)$ and thus, $h({\uparrow}(v^{*}\vee
v))\leqslant k-1$. By inductive hypothesis we have that ${\rm
Ind}({\uparrow}(v^{*}\vee v))< k-1$, proving that ${\rm
Ind}(L)\leqslant k-1$ or equivalently ${\rm Ind}(L)<k$. $\Box$

\section{Sum and product properties for Ind}\label{sec5}

In this section we study properties of the large inductive
dimension. Especially, we study the large inductive dimension of the
linear sum and kinds of products of finite lattices.

\begin{definition}{\rm The \emph{linear sum} of two
posets $(L_1,\leqslant _1)$ and $(L_2,\leqslant _2)$ such that $L_1
\cap L_2=\emptyset$, denoted by $L_1 \oplus L_2$, is the poset $(L_1
\cup L_2,\leqslant)$, where the relation $\leqslant$ is defined as
follows:
\[x\leqslant y \Leftrightarrow \begin{cases} x,y \in L_1 \
\mathrm{and} \ x\leqslant _{1} y\\ x,y \in L_2 \  \mathrm{and} \
x\leqslant _{2} y\\ x\in L_1, \ y\in L_2. \end{cases}\]}
\end{definition}

\noindent Clearly, if $(L_1,\leqslant _1)$ and $(L_2,\leqslant _2)$
are lattices, then $(L_1\oplus L_2,\leqslant)$ is also a lattice.

\begin{remark}{\rm \label{rr1} (1) The lattices $(L_1,\leqslant_1)$ and $(L_3,\leqslant_3)$ of Figures \ref{fig:1} and \ref{fig:9}, respectively, verifies the claim
$\mathrm{Ind} (L_1 \oplus L_3)\neq\mathrm{Ind}(L_3 \oplus L_1)$.
\begin{figure}[H]
\tikzstyle arrowstyle=[scale=1] \tikzstyle
directed=[postaction={decorate,decoration={markings,mark=at position
.6 with {\arrow[arrowstyle]{stealth}}}}]
\begin{center}
\begin{tikzpicture}
\draw[directed] (0,0)--(0,1); \draw[fill] (0,0) circle [radius=0.05]
node [below] {$0_{L_3}$}; \draw[fill] (0,1) circle [radius=0.05]
node [above] {$1_{L_3}$};
\end{tikzpicture}
\end{center}
\caption{The finite lattice $(L_3,\leqslant _3)$} \label{fig:9}
\end{figure}
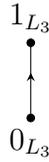
\begin{figure}[H]
\tikzstyle arrowstyle=[scale=1] \tikzstyle
directed=[postaction={decorate,decoration={markings,mark=at position
.6 with {\arrow[arrowstyle]{stealth}}}}]
\begin{center}
\begin{tikzpicture}
\draw[directed] (0,0)--(-1,1); \draw[directed] (0,0)--(1,1);
\draw[directed] (-1,1)--(0,2); \draw[directed] (1,1)--(0,2);
\draw[directed] (0,2)--(0,3); \draw[directed] (0,3)--(0,4);
\draw[directed] (4,0)--(4,1); \draw[directed] (4,1)--(4,2);
\draw[directed] (4,2)--(3,3); \draw[directed] (4,2)--(5,3);
\draw[directed] (3,3)--(4,4); \draw[directed] (5,3)--(4,4);
\draw[fill] (0,0) circle [radius=0.05] node [below] {$0_{L_1}$};
\draw[fill] (-1,1) circle [radius=0.05] node [left] {$x_1$};
\draw[fill] (1,1) circle [radius=0.05] node [right] {$x_2$};
\draw[fill] (0,2) circle [radius=0.05] node [right] {$1_{L_1}$};
\draw[fill] (0,3) circle [radius=0.05] node [right] {$0_{L_3}$};
\draw[fill] (0,4) circle [radius=0.05] node [above] {$1_{L_3}$};
\draw[fill] (4,0) circle [radius=0.05] node [below] {$0_{L_3}$};
\draw[fill] (4,1) circle [radius=0.05] node [right] {$1_{L_3}$};
\draw[fill] (4,2) circle [radius=0.05] node [right] {$0_{L_1}$};
\draw[fill] (3,3) circle [radius=0.05] node [left] {$x_1$};
\draw[fill] (5,3) circle [radius=0.05] node [right] {$x_2$};
\draw[fill] (4,4) circle [radius=0.05] node [above] {$1_{L_1}$};
\end{tikzpicture}
\end{center}
\caption{The finite lattices  $(L_1 \oplus L_3,\leqslant)$ and $(L_3
\oplus L_1,\leqslant ^{*})$} \label{fig:10}
\end{figure}
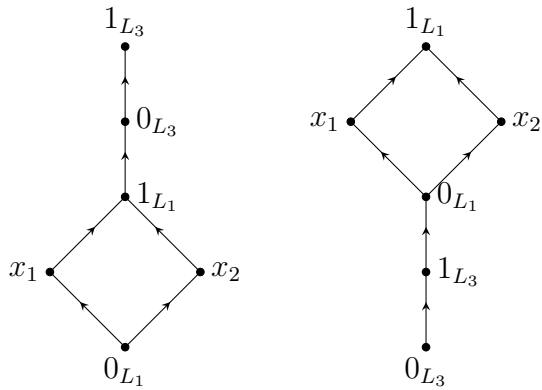
\noindent We have $\mathrm{Ind} (L_1)=\mathrm{Ind} (L_3)=0$,
$\mathrm{Ind} (L_1 \oplus L_3)=0$ and $\mathrm{Ind} (L_3\oplus
L_1)=1$.

\medskip
\noindent (2) The lattice $(L_3\oplus L_1,\leqslant^{*})$ of Figure
\ref{fig:10} also verifies that the assertion
\[\mathrm{Ind}(L_3 \oplus L_1)\leqslant\mathrm{Ind}
(L_3)+\mathrm{Ind}(L_1)\] does not hold for all finite lattices
$L_3$ and $L_1$.}
\end{remark}

\begin{proposition}Let $(L_1,\leqslant_1)$ and $(L_2,\leqslant_2)$ be two finite
lattices such that the pseudocomplement of every element of
$L_2\setminus \{0_{L_{2}}\}$ is $0_{L_{2}}$. Then
\[\mathrm{Ind}(L_2)\leqslant\mathrm{Ind}(L_1\oplus L_2).\]
\end{proposition}
{\bf Proof.} Let $\mathrm{Ind}(L_1\oplus L_2)=k$, where
$k\in\mathbb{N}$, $a\in L_2$ and $v\in L_2$ such that $a\vee
v=1_{L_{2}}$.

(1) If $a=1_{L_{2}}$ and $v\in L_2$, then there exists the element
$u=0_{L_{2}}$ such that $u\leqslant_2 v$ and
$\mathrm{Ind}({\uparrow}(u\vee u^{*}_{L_{2}}))=-1\leqslant k-1$,
where $u^{*}_{L_{2}}$ denotes the pseudocomplement of $u$ in the
lattice $L_{2}$.

(2) If $a=0_{L_{2}}$, then the only element $v\in L_2$ for which
$a\vee v=1_{L_{2}}$ is the element $1_{L_{2}}$. Thus, for this
element $u=1_{L_{2}}$ we have $\mathrm{Ind}({\uparrow}(u\vee
u^{*}_{L_{2}}))=-1\leqslant k-1$.

(3) Let $a\in L_{2}\setminus \{1_{L_{2}},0_{L_{2}}\}$. Then $a,v\in
L_1\oplus L_2$ with $a\vee v=1_{L_1\oplus L_2}$. Since
$\mathrm{Ind}(L_1\oplus L_2)=k$, there exists $u\in L_1\oplus L_2$
such that $u\leqslant v$, $a\vee u=1_{L_1\oplus L_2}$ (and
therefore, $a\vee u=1_{{L}_2}$) and $\mathrm{Ind}({\uparrow}(u\vee
u^{*}_{L_{1}\oplus L_{2}}))\leqslant k-1$, where $u^{*}_{L_{1}\oplus
L_{2}}$ denotes the pseudocomplement of $u$ in the lattice
$L_{1}\oplus L_{2}$.

By the definition of the linear sum $L_{1}\oplus L_{2}$, in order to
have $a\vee u=1_{L_2}$, $u\in L_2$ and thus,
\begin{align*} u^{*}_{L_{1}\oplus L_{2}}&=\bigvee\{x\in L_{1}\oplus L_{2}:x\wedge u=0_{L_{1}\oplus L_{2}}\}\\
&=\bigvee\{x\in L_{1}\oplus L_{2}:x\wedge u=0_{L_{1}}\}\\
&=0_{L_{1}}.
\end{align*}
Therefore, \[\mathrm{Ind}({\uparrow}(u\vee u^{*}_{L_{1}\oplus
L_{2}}))=\mathrm{Ind}({\uparrow}(u\vee
0_{L_{1}}))=\mathrm{Ind}({\uparrow}u)\leqslant k-1\] and thus,
\[\mathrm{Ind}({\uparrow}(u\vee
u^{*}_{L_{2}}))=\mathrm{Ind}({\uparrow}(u\vee
0_{L_{2}}))=\mathrm{Ind}({\uparrow}u)\leqslant k-1.\] Hence,
$\mathrm{Ind}(L_2)\leqslant k$. $\Box$

\begin{definition}{\rm If $(L_1,\leqslant _1)$ and $(L_2,\leqslant
_2)$ are two posets, then their \emph{Cartesian product} is the
poset $(L_1\times L_2,\leqslant)$, where the relation $\leqslant$ is
defined as follows: $(x_1,y_1)\leqslant (x_2,y_2)\Leftrightarrow x_1
\leqslant _1 x_2 \ \mathrm{and} \ y_1 \leqslant _2 y_2$.}
\end{definition}

If $(L_1,\leqslant _1)$ and $(L_2,\leqslant _2)$ are finite
lattices, then $(L_1 \times L_2,\leqslant)$ is also a finite
lattice. Also, if $(x_1,y_1),(x_2,y_2)\in L_1 \times L_2$, then:\\
(1) $(x_1,y_1)\vee (x_2,y_2)=(x_1\vee x_2,y_1\vee y_2)$ and\\
(2) $(x_1,y_1)\wedge (x_2,y_2)=(x_1\wedge x_2,y_1\wedge y_2)$.

\begin{lemma} \cite{CM,GMPS} \label{lem1} If $L_1$ and $L_2$ are finite lattices and
$(x,y)\in L_1\times L_2$, then the following are satisfied:
\begin{enumerate}
\item ${\uparrow}(x,y)={\uparrow}x \times {\uparrow}y,$
\item $(x,y)^{*}=(x^{*},y^{*})$.
\end{enumerate}
\end{lemma}

\begin{theorem} \label{mtheo2} Let $(L_1,\leqslant_1)$ and $(L_2,\leqslant_2)$ be finite lattices.
Then,
\[\mathrm{Ind}(L_1 \times L_2)=\max\{\mathrm{Ind}(L_1),\mathrm{Ind}(L_2)\}.\]
\end{theorem}
{\bf Proof.} In order to prove the inequality
\[\mathrm{Ind}(L_1 \times L_2)\leqslant \max\{\mathrm{Ind}(L_1), \mathrm{Ind}(L_2)\},\]
we shall apply induction with respect to the number
\[k(L_1,L_2)=\max\{\mathrm{Ind}(L_1), \mathrm{Ind}(L_2)\}.\]
If $k(L_1,L_2)=-1$, then our inequality holds. Moreover, if
$\mathrm{Ind}(L_1)=-1$ and $\mathrm{Ind}(L_2)\geqslant 0$
(respectively, $\mathrm{Ind}(L_2)=-1$ and
$\mathrm{Ind}(L_1)\geqslant 0$), then the lattices $L_1 \times L_2$
and $L_2$ (respectively, $L_1 \times L_2$ and $L_1$) are isomorphic
and therefore, $\mathrm{Ind}(L_1 \times L_2)=\mathrm{Ind}(L_2)$
(respectively, $\mathrm{Ind}(L_1 \times L_2)=\mathrm{Ind}(L_1)$).

We assume that the inequality holds for every pair of finite
lattices $L_1$ and $L_2$ with $k(L_1,L_2)<k$, where $k\geqslant 0$,
and we consider finite lattices $L_1$ and $L_2$ such that
$\mathrm{Ind}(L_1)=n\geqslant 0$, $\mathrm{Ind}(L_2)=m\geqslant 0$
and $\max\{n,m\}=k$. We shall prove that $\mathrm{Ind}(L_1 \times
L_2)\leqslant k$.

Let $(a,b)\in L_1\times L_2$ and $(u,v)\in L_1\times L_2$ such that
$(a,b)\vee (u,v)=(1_{L_{1}},1_{L_{2}})$. Then $a,u\in L_1$ with
$a\vee u=1_{L_{1}}$ and $b,v\in L_2$ with $b\vee v=1_{L_{2}}$. Since
$\mathrm{Ind}(L_1)=n$, there exists $x\in L_1$ with $x\leqslant_1
u$, $x\vee a=1_{L_{1}}$ and \[\mathrm{Ind}({\uparrow}(x\vee
x^{*}))\leqslant n-1.\] Similarly, since $\mathrm{Ind}(L_2)=m$,
there exists $y\in L_2$ with $y\leqslant_2 v$, $y\vee b=1_{L_{2}}$
and
\[\mathrm{Ind}({\uparrow}(y\vee y^{*}))\leqslant m-1.\]

We consider the element $(x,y)\in L_1\times L_2$. Then
$(x,y)\leqslant (u,v)$ and $(a,b)\vee (x,y)=(1_{L_{1}},1_{L_{2}})$.
We shall prove that \[\mathrm{Ind}({\uparrow}((x,y)^{*}\vee
(x,y)))\leqslant \max\{n,m\}-1.\] Considering Lemma \ref{lem1} we
have that
\begin{align*} \mathrm{Ind}({\uparrow}((x,y)^{*}\vee (x,y)))&=\mathrm{Ind}({\uparrow}((x^{*},y^{*})\vee (x,y)))\\
&=\mathrm{Ind}({\uparrow}(x^{*}\vee x,y^{*}\vee y))\\
&=\mathrm{Ind}({\uparrow}(x^{*}\vee x) \times {\uparrow}(y^{*}\vee
y)).
\end{align*}
By inductive hypothesis, we have that
\begin{align*}\mathrm{Ind}({\uparrow}(x^{*}\vee x) \times {\uparrow}(y^{*}\vee
y))&\leqslant \max\{\mathrm{Ind}({\uparrow}(x^{*}\vee x)),
\mathrm{Ind}({\uparrow}(y^{*}\vee y))\}\\&\leqslant\max\{n-1,m-1\}\\
&=\max\{n,m\}-1.\nonumber
\end{align*}
Therefore, $\mathrm{Ind}(L_1 \times L_2)\leqslant \max\{n,m\}$.

We shall prove the inequality
\[\max\{\mathrm{Ind}(L_1),\mathrm{Ind}(L_2)\}\leqslant\mathrm{Ind}(L_1
\times L_2).\] If $\mathrm{Ind}(L_1 \times L_2)=-1$, then $L_1\times
L_2=\{1_{L_1\times L_2}\}$ and thus, $L_1=\{1_{L_1}\}$ and
$L_2=\{1_{L_2}\}$ and the inequality holds.

We suppose that the inequality holds for every pair of finite
lattices $L_1$ and $L_2$ with $\mathrm{Ind}(L_1 \times L_2)<k$ and
we shall prove it for $k$. Let $\mathrm{Ind}(L_1 \times L_2)=k$. We
shall prove that $\mathrm{Ind}(L_1)\leqslant k$ and
$\mathrm{Ind}(L_2)\leqslant k$.

Let $a\in L_1$ and $u\in L_1$ with $a\vee u=1_{L_{1}}$. Let also
$b\in L_2$ and $v\in L_2$ with $b\vee v=1_{L_{2}}$. Then $(a,b),
(u,v)\in L_1\times L_2$ with $(a,b)\vee
(u,v)=(1_{L_{1}},1_{L_{2}})$. Since $\mathrm{Ind}(L_1 \times
L_2)=k$, there exists $(x,y)\in L_1\times L_2$ such that
$(x,y)\leqslant (u,v)$, $(a,b)\vee (x,y)=(1_{L_{1}},1_{L_{2}})$ and
\[\mathrm{Ind}({\uparrow}((x,y)^{*}\vee (x,y)))\leqslant k-1.\] That
is, $x\in L_1$ with $x\leqslant_1 u$ and $a\vee x=1_{L_{1}}$ and
$y\in L_2$ with $y\leqslant_2 v$ and $b\vee y=1_{L_{2}}$. It
suffices to prove that $\mathrm{Ind}({\uparrow}(x^{*}\vee
x))\leqslant k-1$ and $\mathrm{Ind}({\uparrow}(y^{*}\vee
y))\leqslant k-1$.

Considering Lemma \ref{lem1} we have that
\begin{align*} \mathrm{Ind}({\uparrow}((x,y)^{*}\vee (x,y)))&=\mathrm{Ind}({\uparrow}((x^{*},y^{*})\vee (x,y)))\\
&=\mathrm{Ind}({\uparrow}(x^{*}\vee x,y^{*}\vee y))\\
&=\mathrm{Ind}({\uparrow}(x^{*}\vee x) \times {\uparrow}(y^{*}\vee
y)).
\end{align*}
By inductive hypothesis,
\[\max\{\mathrm{Ind}({\uparrow}(x^{*}\vee x)), \mathrm{Ind}({\uparrow}(y^{*}\vee
y))\}\leqslant\mathrm{Ind}({\uparrow}(x^{*}\vee x) \times
{\uparrow}(y^{*}\vee y)) \leqslant k-1.\] Therefore,
$\mathrm{Ind}({\uparrow}(x^{*}\vee x))\leqslant k-1$ and $
\mathrm{Ind}({\uparrow}(y^{*}\vee y))\leqslant k-1$. Thus,
$\mathrm{Ind}(L_1)\leqslant k$ and $\mathrm{Ind}(L_2)\leqslant k$.
$\Box$

\medskip
Thus, in contrast to Remark \ref{rr1} for the linear sum, Theorem
\ref{mtheo2} tends to the following corollaries which show the
respect of Ind to a kind of Cartesian product theorem and the
property of commutativity.

\begin{corollary} Let $(L_1,\leqslant_1)$ and $(L_2,\leqslant_2)$ be finite lattices.
Then,
\begin{enumerate}
\item
$\mathrm{Ind}(L_1 \times L_2)\leqslant
\mathrm{Ind}(L_1)+\mathrm{Ind}(L_2)$,
\item
$\mathrm{Ind}(L_1 \times L_2)=\mathrm{Ind}(L_2 \times L_1)$.
\end{enumerate}
\end{corollary}

\begin{definition}{\rm For two lattices
$(L_1,\leqslant_1)$ and $(L_2,\leqslant _2)$ the \emph{lexicographic
product} $L_1 \diamond L_2$ is the lattice $(L_1\times
L_2,\leqslant)$, where the relation $\leqslant$ is defined as
follows: \[(x_1,y_1)\leqslant (x_2,y_2) \Leftrightarrow
\begin{cases} x_1 < _{1} x_2 \ \mathrm{or}\\ x_1 =x_2 \ \mathrm{and}
\ y_1 \leqslant _{2} y_2. \end{cases}\]}
\end{definition}

\noindent By Definition 4.8 of \cite{DGMM} for every $(x,y),
(x^{\prime}, y^{\prime})\in L_1\diamond L_2$ we have that
\[(x,y)\wedge (x^{\prime}, y^{\prime})=\begin{cases}(x,y), \
\mathrm{if} \
x<x^{\prime}\\ (x^{\prime},y^{\prime}), \ \mathrm{if} \ x^{\prime}<x\\
(x, y\wedge y^{\prime}), \ \mathrm{if} \ x=x^{\prime}\\ (x\wedge
x^{\prime}, 1_{L_{2}}), \ \mathrm{if} \ x\parallel
x^{\prime}\end{cases}\] and \[(x,y)\vee(x^{\prime},
y^{\prime})=\begin{cases}(x,y), \ \mathrm{if} \
x^{\prime}<x\\ (x^{\prime},y^{\prime}), \ \mathrm{if} \ x<x^{\prime}\\
(x, y\vee y^{\prime}), \ \mathrm{if} \ x=x^{\prime}\\ (x\vee
x^{\prime}, 0_{L_{2}}), \ \mathrm{if} \ x\parallel
x^{\prime}.\end{cases}\]

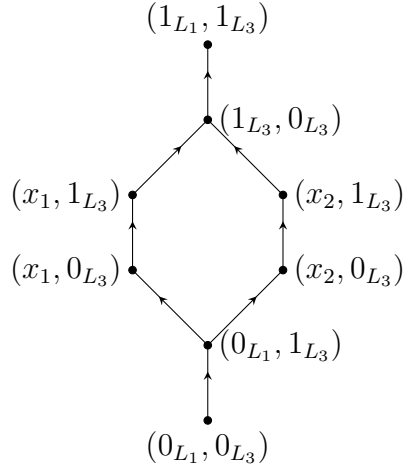
\begin{figure}[H]
\tikzstyle arrowstyle=[scale=1] \tikzstyle
directed=[postaction={decorate,decoration={markings,mark=at position
.65 with {\arrow[arrowstyle]{stealth}}}}]
\begin{center}
\begin{tikzpicture}
\draw[directed] (0,-1)--(0,0); \draw[directed] (0,0)--(-1,1);
\draw[directed] (0,0)--(1,1); \draw[directed] (-1,1)--(-1,2);
\draw[directed] (1,1)--(1,2); \draw[directed] (-1,2)--(0,3);
\draw[directed] (1,2)--(0,3); \draw[directed] (0,3)--(0,4);
\draw[fill] (0,-1) circle [radius=0.05] node [below]
{$(0_{L_1},0_{L_3})$}; \draw[fill] (0,0) circle [radius=0.05] node
[right] {$(0_{L_1},1_{L_3})$}; \draw[fill] (-1,1) circle
[radius=0.05] node [left] {$(x_1,0_{L_3})$}; \draw[fill] (1,1)
circle [radius=0.05] node [right] {$(x_2,0_{L_3})$}; \draw[fill]
(-1,2) circle [radius=0.05] node [left] {$(x_1,1_{L_{3}})$};
\draw[fill] (1,2) circle [radius=0.05] node [right]
{$(x_2,1_{L_{3}})$}; \draw[fill] (0,3) circle [radius=0.05] node
[right] {$(1_{L_3},0_{L_{3}})$}; \draw[fill] (0,4) circle
[radius=0.05] node [above] {$(1_{L_1},1_{L_3})$};
\end{tikzpicture}
\end{center}
\caption{The finite lattice $(L_1\diamond L_3,\leqslant)$}
\label{fig:23}
\end{figure}

\begin{remark}{\rm (1) The lattices $(L_1,\leqslant_1)$ and $(L_3,\leqslant_3)$ of Figures \ref{fig:1} and \ref{fig:9}, respectively, verifies that the
inequality
\[\mathrm{Ind}(L_3 \diamond L_1)\leqslant\mathrm{Ind}(L_3)+\mathrm{Ind}(L_1)\] does not hold for all finite lattices $L_1$
and $L_3$.

We have that $\mathrm{Ind}(L_1)=\mathrm{Ind}(L_3)=0$. For the finite
lattice $(L_3\diamond L_1,\leqslant)$, which is represented by the
diagram of Figure \ref{fig:11}, we have that
$\mathrm{Ind}(L_3\diamond L_1)=1$.
\begin{figure}[H]
\tikzstyle arrowstyle=[scale=1] \tikzstyle
directed=[postaction={decorate,decoration={markings,mark=at position
.65 with {\arrow[arrowstyle]{stealth}}}}]
\begin{center}
\begin{tikzpicture}
\draw[directed] (0,0)--(-1,1); \draw[directed] (0,0)--(1,1);
\draw[directed] (-1,1)--(0,2); \draw[directed] (1,1)--(0,2);
\draw[directed] (0,2)--(0,3); \draw[directed] (0,3)--(-1,4);
\draw[directed] (0,3)--(1,4); \draw[directed] (-1,4)--(0,5);
\draw[directed] (1,4)--(0,5); \draw[fill] (0,0) circle [radius=0.05]
node [below] {$(0_{L_3},0_{L_1})$}; \draw[fill] (-1,1) circle
[radius=0.05] node [left] {$(0_{L_3},x_1)$}; \draw[fill] (1,1)
circle [radius=0.05] node [right] {$(0_{L_3},x_2)$}; \draw[fill]
(0,2) circle [radius=0.05] node [right] {$(0_{L_3},1_{L_{1}})$};
\draw[fill] (0,3) circle [radius=0.05] node [right]
{$(1_{L_3},0_{L_{1}})$}; \draw[fill] (-1,4) circle [radius=0.05]
node [left] {$(1_{L_3},x_1)$}; \draw[fill] (1,4) circle
[radius=0.05] node [right] {$(1_{L_3},x_2)$}; \draw[fill] (0,5)
circle [radius=0.05] node [above] {$(1_{L_3},1_{L_1})$};
\end{tikzpicture}
\end{center}
\caption{The finite lattice $(L_3\diamond L_1,\leqslant)$}
\label{fig:11}
\end{figure}
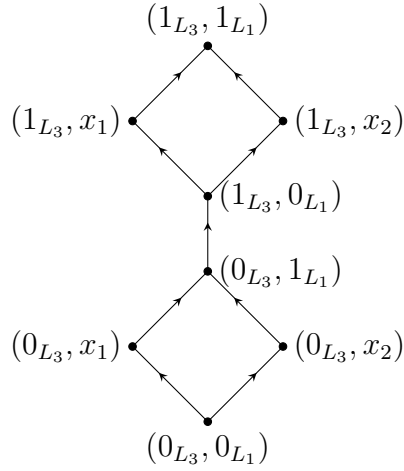

\noindent (2) The lattices $(L_1,\leqslant_1)$ and
$(L_3,\leqslant_3)$ of Figures \ref{fig:1} and \ref{fig:9},
respectively, verifies that the inequality
\[\mathrm{Ind}(L_3 \diamond L_1)=\mathrm{Ind}(L_1\diamond L_3)\] does not hold for all finite lattices $L_1$
and $L_3$.

As we have seen in Figure \ref{fig:11}, $\mathrm{Ind}(L_3\diamond
L_1)=1$. However, $\mathrm{Ind}(L_1\diamond L_3)=0$. Indeed, the
diagram of $L_1\diamond L_3$ is given in Figure \ref{fig:23}, from
which we can easily see that $\mathrm{Ind}(L_1\diamond L_3)=0$.}
\end{remark}

\begin{proposition} Let $(L_1,\leqslant_1)$ and $(L_2,\leqslant_2)$ be two finite lattices such that the pseudocomplement of every element of
$L_2\setminus\{0_{L_2}\}$ is $0_{L_2}$. Then
\[\mathrm{Ind}(L_2)\leqslant \mathrm{Ind}(L_1\diamond L_2).\]
\end{proposition}
{\bf Proof.} Let $\mathrm{Ind}(L_1\diamond L_2)=k$,
$k\in\mathbb{N}$, $a\in L_2$ and $v\in L_2$ with $a\vee
v=1_{L_{2}}$.

(1) If $a=1_{L_{2}}$, then for every $v\in L_2$ we have $a\vee
v=1_{L_{2}}$ and thus, there exists the element
$u^{\prime}=0_{L_{2}}$ such that $u^{\prime}\leqslant_2 v$, $a\vee
u^{\prime}=1_{L_{2}}$ and \[\mathrm{Ind}({\uparrow}(u^{\prime}\vee
(u^{\prime})^{*}_{L_{2}}))=-1\leqslant k-1,\] where
$(u^{\prime})^{*}_{L_{2}}$ denotes the pseudocomplement of
$u^{\prime}$ in $L_2$.

(2) If $a=0_{L_{2}}$, then for the only element $v=1_{L_2}$ we have
$a\vee v=1_{L_{2}}$. Thus, there exists the element
$u^{\prime}=1_{L_{2}}$ such that $u^{\prime}\leqslant_2 v$, $a\vee
u^{\prime}=1_{L_{2}}$ and \[\mathrm{Ind}({\uparrow}(u^{\prime}\vee
(u^{\prime})^{*}_{L_{2}}))=-1\leqslant k-1.\]

(3) Let $a\in L_2\setminus \{0_{L_{2}},1_{L_{2}}\}$. Then
$(1_{L_{1}},a)$ and $(1_{L_{1}},v)$ belong to $L_1\diamond L_2$ and
$(1_{L_{1}},a)\vee (1_{L_{1}},v)=(1_{L_{1}},1_{L_{2}})$. Since
$\mathrm{Ind}(L_1\diamond L_2)=k$, there exists $(u,u^{\prime})\in
L_1\diamond L_2$ such that $(u,u^{\prime})\leqslant (1_{L_{1}},v)$,
$(u,u^{\prime})\vee (1_{L_{1}},a)=(1_{L_{1}},1_{L_{2}})$ and
\[\mathrm{Ind}({\uparrow}((u,u^{\prime})\vee
(u,u^{\prime})^{*}_{L_{1}\diamond L_2}))\leqslant k-1,\] where
$(u,u^{\prime})^{*}_{L_{1}\diamond L_{2}}$ denotes the
pseudocomplement of $(u,u^{\prime})$ in $L_{1}\diamond L_{2}$. Since
$(u,u^{\prime})\leqslant (1_{L_{1}},v)$, $a\neq 1_{L_{2}}$ and
$(u,u^{\prime})\vee (1_{L_{1}},a)=(1_{L_{1}},1_{L_{2}})$, we have
that $u=1_{L_{1}}$, $u^{\prime}\leqslant_{2} v$ and $u^{\prime}\vee
a=1_{L_{2}}$. Since $(u')_{L_2}^* = 0_{L_2}$ by the assumption of
the proposition, we have $(1_{L_{1}},u^{\prime})^{*}_{L_{1}\diamond
L_{2}}=(0_{L_{1}},0_{L_{2}})$ and thus,
\[\mathrm{Ind}({\uparrow}(u,u^{\prime}))=\mathrm{Ind}({\uparrow}((u,u^{\prime})\vee (u,u^{\prime})^{*}_{L_{1}\diamond
L_{2}}))\leqslant k-1.\] Since
${\uparrow}(1_{L_1},u^{\prime})={\uparrow}1_{L_1}\times{\uparrow}u^{\prime}$,
we have that
\[\mathrm{Ind}({\uparrow}u^{\prime})=\mathrm{Ind}(\{1_{L_1}\}\times
{\uparrow}u^{\prime})=\mathrm{Ind}({\uparrow}(1_{L_1},u^{\prime}))\leqslant
k-1\] and thus,
\[\mathrm{Ind}({\uparrow}(u^{\prime}\vee
(u^{\prime})^{*}_{L_{2}}))=\mathrm{Ind}({\uparrow}(u^{\prime}\vee
0_{L_{2}}))=\mathrm{Ind}({\uparrow}u^{\prime})\leqslant k-1.\]
Therefore, $\mathrm{Ind}(L_2)\leqslant k$. $\Box$

\begin{figure}[H]
\tikzstyle arrowstyle=[scale=1] \tikzstyle
directed=[postaction={decorate,decoration={markings,mark=at position
.6 with {\arrow[arrowstyle]{stealth}}}}]
\begin{center}
\begin{tikzpicture}
\draw [directed] (0,0) -- (0,1); \draw [directed] (0,1) -- (0,2);
\draw [directed] (0,2) -- (0,3); \draw [directed] (0,3) -- (-1,4);
\draw [directed] (0,3) -- (1,4); \draw [directed] (-1,4) -- (-1,5);
\draw [directed] (1,4) -- (1,5); \draw [directed] (-1,5) -- (0,6);
\draw [directed] (1,5) -- (0,6); \draw [directed] (0,6) -- (0,7);
\draw[fill] (0,0) circle [radius=0.05] node [below]
{$(0_{L_{2}},0_{L_{3}})$}; \draw[fill] (0,1) circle [radius=0.05]
node [right] {$(0_{L_{2}},1_{L_{3}})$}; \draw[fill] (0,2) circle
[radius=0.05] node [right] {$(y_1,0_{L_{3}})$}; \draw[fill] (0,3)
circle [radius=0.05] node [right] {$(y_1,1_{L_{3}})$}; \draw[fill]
(-1,4) circle [radius=0.05] node [left] {$(y_2,0_{L_{3}})$};
\draw[fill] (1,4) circle [radius=0.05] node [right]
{$(y_3,0_{L_3})$}; \draw[fill] (-1,5) circle [radius=0.05] node
[left] {$(y_2,1_{L_{3}})$}; \draw[fill] (1,5) circle [radius=0.05]
node [right] {$(y_3,1_{L_3})$}; \draw[fill] (0,6) circle
[radius=0.05] node [right] {$(1_{L_{2}},0_{L_{3}})$}; \draw[fill]
(0,7) circle [radius=0.05] node [right] {$(1_{L_{2}},1_{L_{3}})$};
\end{tikzpicture}
\end{center}
\caption{The finite lattice $(L_2\diamond L_3,\leqslant)$}
\label{fig:12}
\end{figure}
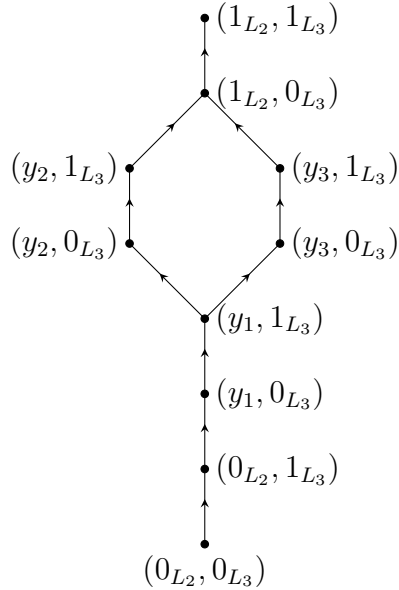

\begin{remark}{\rm The dimension Ind of the Cartesian
product $L_1\times L_2$ of two finite lattices $L_1$ and $L_2$
differs from the dimension Ind of their lexicographic product.

\medskip
\noindent (1) For the lattices $(L_1,\leqslant_1)$ and
$(L_3,\leqslant_3)$ of Figures \ref{fig:1} and \ref{fig:9},
respectively, we have $\mathrm{Ind}(L_1)=\mathrm{Ind}(L_3)=0$ and by
Theorem \ref{mtheo2}, $\mathrm{Ind}(L_3\times L_1)=0$. For their
lexicographic product (see Figure \ref{fig:11}) we have
$\mathrm{Ind}(L_3\diamond L_1)=1$.

\medskip
\noindent (2) For the lattices $(L_2,\leqslant_2)$ and
$(L_3,\leqslant_3)$ of Figures \ref{fig:1} and \ref{fig:9},
respectively, we have $\mathrm{Ind}(L_2)=1$, $\mathrm{Ind}(L_3)=0$
and by Theorem \ref{mtheo2}, $\mathrm{Ind}(L_2\times L_3)=1$. For
their lexicographic product (see Figure \ref{fig:12}) we have
$\mathrm{Ind}(L_2\diamond L_3)=0$.}
\end{remark}

\begin{definition}{\rm \cite{Bae1,Ben} The \emph{rectangular product} of two
lattices $(L_1,\leqslant _1)$ and $(L_2,\leqslant _2)$ is the
lattice $(L_1 \Box L_2,\leqslant)$, where
\[L_1 \Box L_2=\{(x,y)\in L_1\times L_2: x\neq 0_{L_1} \ \mathrm{and} \ y\neq
0_{L_2}\} \cup \{(0_{L_1},0_{L_2})\}\] and the relation $\leqslant$
is defined as follows:
\[(x_1,y_1)\leqslant (x_2,y_2)\Leftrightarrow x_1
\leqslant _{1} x_2 \ \mathrm{and} \ y_1 \leqslant _2 y_2.\]}
\end{definition}

That is, $L_1 \Box L_2$ is the subset of $L_1 \times L_2$ obtained
by removing the points of the form $(x,0_{L_2})$ and $(0_{L_1},y)$
with $x\neq 0_{L_1}$ and $y\neq 0_{L_2}$. Also, if
$(x_1,y_1),(x_2,y_2)\in L_1 \Box L_2$, then:\\
(1) $(x_1,y_1)\vee (x_2,y_2)=(x_1\vee x_2,y_1\vee y_2)$ and\\
(2) $(x_1,y_1)\wedge (x_2,y_2)=\begin{cases} (x_1\wedge
x_2,y_1\wedge y_2), \ \mbox{if} \ x_1\wedge x_2\neq 0_{L_1} \
\mbox{and} \ y_1\wedge y_2\neq 0_{L_2}\\
(0_{L_1},0_{L_2}),\ \mbox{otherwise}.\end{cases}$

\medskip
In order to study the relation of the rectangular product of two
finite lattices with their large inductive dimension, we recall the
following fact.

\begin{proposition} \cite{GMPS} \label{pseudo} Let $L_1$ and $L_2$ be two finite lattices and
$(a,b)\in L_1 \Box L_2$. Then
\[(a,b)^{*}=\begin{cases} (a^{*},b^{*})=(0_{L_{1}},0_{L_{2}}), & \ \mathrm{if} \ a^{*}=0_{L_{1}} \ \mathrm{and} \ b^{*}=0_{L_{2}}\\
(1_{L_{1}},1_{L_{2}}), & \ \mathrm{if} \ a^{*}\neq 0_{L_{1}} \ \mathrm{and} \ b^{*}\neq 0_{L_{2}}\\
(1_{L_{1}},b^{*}), & \ \mathrm{if} \ a^{*}=0_{L_{1}} \ \mathrm{and} \ b^{*}\neq 0_{L_{2}}\\
(a^{*},1_{L_{2}}), & \ \mathrm{if} \ a^{*}\neq 0_{L_{1}} \
\mathrm{and} \ b^{*}=0_{L_{2}}.
\end{cases}\]
\end{proposition}

Up to now, the above discussions related the dimension ${\rm Ind}$
and the properties of linear sum, Cartesian product and
lexicographic product tend to study the behavior of this dimension
without mentioning additional lattice properties, verifying that Ind
respects such properties. However, this seems to fail for the
rectangular product.

Let us observe that in Figure \ref{fig:3} we have seen a finite
lattice $L$ for which
$\mathrm{Ind}({\uparrow}x_2)\nleqslant\mathrm{Ind}(L)$. This fact
leads us to insert a new lattice property as follows. This property
will be useful in order to succeed a kind of rectangular product
theorem given in Theorem \ref{mtheo3}.

\begin{definition}{\rm Let $L$ be a finite lattice. If the sublattice
property of the form
$\mathrm{Ind}({\uparrow}x)\leqslant\mathrm{Ind}(L)$ holds for every
$x\in L$, then we say that $L$ has the {\it SP-property}.}
\end{definition}

\begin{theorem} \label{mtheo3} Let $(L_1,\leqslant_1)$ and $(L_2,\leqslant_2)$ be two finite
lattices such that the SP-property holds for all ${\uparrow}\kappa$
and ${\uparrow}\lambda$, where $\kappa\in L_1$ and $\lambda\in L_2$.
Then \[\mathrm{Ind}(L_1\Box L_2)\leqslant
\max\{\mathrm{Ind}(L_1),\mathrm{Ind}(L_2)\}+1.\]
\end{theorem}
{\bf Proof.} Let $\mathrm{Ind}(L_1)=n$ and $\mathrm{Ind}(L_2)=m$,
where $n,m\in\mathbb{N}$. We shall prove that $\mathrm{Ind}(L_1\Box
L_2)\leqslant \max\{n,m\}+1$. Let $(a,b)\in L_1\Box L_2$. Then we
study the following cases:

\noindent (1) If $(a,b)=(0_{L_{1}},0_{L_{2}})$, then for every
$(u,v)\in L_1\Box L_2$ with $(a,b)\vee (u,v)=(1_{L_{1}},1_{L_{2}})$
(that is, $(u,v)=(1_{L_{1}},1_{L_{2}})$), there is
$(x,y)=(1_{L_{1}},1_{L_{2}})\in L_1\Box L_2$ such that
$(x,y)\leqslant (u,v)$, $(a,b)\vee (x,y)=(1_{L_{1}},1_{L_{2}})$ and
\[\mathrm{Ind}({\uparrow}((x,y)\vee (x,y)^{*}))=-1<\max\{n,m\}.\]

\noindent (2) Let $a\neq 0_{L_{1}}$, $b\neq 0_{L_{2}}$.

(2a) If $(a,b)=(1_{L_{1}},1_{L_{2}})$ and $(u,v)\in L_1\Box L_2$
such that $(a,b)\vee (u,v)=(1_{L_{1}},1_{L_{2}})$, then there exists
$(x,y)=(0_{L_{1}},0_{L_{2}})\in L_1\Box L_2$ such that
$(x,y)\leqslant (u,v)$, $(a,b)\vee (x,y)=(1_{L_{1}},1_{L_{2}})$ and
\[\mathrm{Ind}({\uparrow}((x,y)\vee (x,y)^{*}))=-1<\max\{n,m\}.\]

(2b) Let $(a,b)\in (L_1\Box L_2)\setminus\{(0_{L_{1}},0_{L_{2}}),
(1_{L_{1}},1_{L_{2}})\}$ and $(u,v)\in L_1\Box L_2$ such that
$(a,b)\vee (u,v)=(1_{L_{1}},1_{L_{2}})$. We shall prove that there
exists $(x,y)\in L_1\Box L_2$ such that $(x,y)\leqslant (u,v)$,
$(a,b)\vee (x,y)=(1_{L_{1}},1_{L_{2}})$ and
\[\mathrm{Ind}({\uparrow}((x,y)\vee
(x,y)^{*}))\leqslant\max\{n,m\}.\] We consider as $(x,y)$ the
element $(u,v)$. Then by Proposition \ref{pseudo} we have the
following cases:

(A) $(u,v)\vee (u,v)^{*}=(u,v)\vee(u^{*},1_{L_{2}})=(u\vee u^{*},
1_{L_{2}})$ and thus, by Lemma \ref{lem1} and Theorem \ref{mtheo2},
we have that \[\mathrm{Ind}({\uparrow}((u,v)\vee
(u,v)^{*}))=\max\{\mathrm{Ind}({\uparrow}(u\vee u^{*})),-1\}.\] By
SP-property for the lattice $L_1$ we also have that
$\mathrm{Ind}({\uparrow}(u\vee u^{*}))\leqslant{\rm Ind}(L_1)$ and
thus, \[\mathrm{Ind}({\uparrow}((u,v)\vee
(u,v)^{*}))\leqslant\max\{n,-1\}.\]

(B) $(u,v)\vee (u,v)^{*}=(u,v)\vee
(1_{L_{1}},v^{*})=(1_{L_{1}},v\vee v^{*})$ and thus, this case is
similar to Case (A).

(C) $(u,v)\vee (u,v)^{*}=(u,v)\vee (0_{L_{1}},0_{L_{2}})=(u,v)$ and
thus, by Lemma \ref{lem1} and Theorem \ref{mtheo2}, we have that
\[\mathrm{Ind}({\uparrow}((u,v)\vee
(u,v)^{*}))=\max\{\mathrm{Ind}({\uparrow}u),\mathrm{Ind}({\uparrow}v)\}.\]
By SP-property for the lattices $L_1$ and $L_2$ we also have
$\mathrm{Ind}({\uparrow}u)\leqslant{\rm Ind}(L_1)$ and
$\mathrm{Ind}({\uparrow}v)\leqslant{\rm Ind}(L_2)$. Thus,
\[\mathrm{Ind}({\uparrow}((u,v)\vee (u,v)^{*}))\leqslant\max\{n,m\}.\]

(D) $(u,v)\vee (u,v)^{*}=(u,v)\vee
(1_{L_{1}},1_{L_{2}})=(1_{L_1},1_{L_{2}})$ and thus, \[\mathrm{
Ind}({\uparrow}((u,v)\vee (u,v)^{*}))=-1<\max\{n,m\}.\]

\noindent Therefore, in each case we have that
\[\mathrm{Ind}({\uparrow}((u,v)\vee
(u,v)^{*}))\leqslant\max\{n,m\}.\] Thus, $$\mathrm{Ind}(L_1\Box
L_2)\leqslant \max\{n,m\}+1.\ \Box$$

\begin{corollary} Let $(L_1,\leqslant_1)$ and $(L_2,\leqslant_2)$ be two finite
lattices such that the SP-property holds for all ${\uparrow}\kappa$
and ${\uparrow}\lambda$, where $\kappa\in L_1$ and $\lambda\in L_2$.
Then \[\mathrm{Ind}(L_1\Box L_2)\leqslant
\mathrm{Ind}(L_1)+\mathrm{Ind}(L_2)+1.\]
\end{corollary}

\begin{theorem} \label{thm10} Let $(L_1,\leqslant_1)$ and $(L_2,\leqslant_2)$ be two finite
lattices such that the pseudocomplement of every $w\in
L_i\setminus\{0_{L_{i}}\}$ is $0_{L_{i}}$, for $i=1,2$. Then
\[\mathrm{Ind}(L_1\times L_2)\leqslant\mathrm{Ind}(L_1\Box L_2).\]
\end{theorem}
{\bf Proof.} By Theorem \ref{mtheo2}, it suffices to prove that
\[\max\{\mathrm{Ind}(L_1),\mathrm{Ind}(L_2)\}\leqslant\mathrm{Ind}(L_1\Box
L_2).\] Let $\mathrm{Ind}(L_1\Box L_2)=k$, where $k\in\mathbb{N}$.
Then we shall prove that $\mathrm{Ind}(L_i)\leqslant k$, for
$i=1,2$.

Firstly, we shall prove that
$\mathrm{Ind}(L_1)\leqslant\mathrm{Ind}(L_1\Box L_2)$. Let $a\in
L_1$ and $v\in L_1$ with $a\vee v=1_{L_{1}}$.

(1) If $a=1_{L_{1}}$ and $v\in L_1$, we consider on the element $u=
0_{L_1}$. Then, $u\leqslant_1 v$, $a \vee u = 1_{L_1}$ and
$\mathrm{Ind}(\uparrow(u \vee u^*)) = \mathrm{Ind}(\{1_{L_1}\}) = -1
\leqslant k-1$.

(2) If $a=0_{L_{1}}$, then the only element $v\in L_1$ for which
$a\vee v=1_{L_{1}}$ is the element $1_{L_{1}}$. Thus, for this
element $u=1_{L_{1}}$ we have $\mathrm{Ind}({\uparrow}(u\vee
u^{*}))=-1\leqslant k-1$.

(3) Let $a\in L_1\setminus\{0_{L_{1}},1_{L_{1}}\}$ and $v\in L_1$
with $a\vee v=1_{L_{1}}$. Then $v\neq 0_{L_{1}}$ and thus,
$(a,1_{L_{2}}), (v,1_{L_{2}})\in L_1\Box L_2$ with
$(a,1_{L_{2}})\vee (v,1_{L_{2}})=(1_{L_{1}},1_{L_{2}})$. Since
$\mathrm{Ind}(L_1\Box L_2)=k$, there exists $(x,y)\in L_1\Box L_2$
such that $(x,y)\leqslant (v,1_{L_{2}})$, $(x,y)\vee
(a,1_{L_{2}})=(1_{L_{1}},1_{L_{2}})$ and
\[\mathrm{Ind}({\uparrow}((x,y)\vee (x,y)^{*}))\leqslant k-1.\] Then
$x\leqslant_1 v$ and $x\vee a=1_{L_{1}}$. Also, by assumption we
have \[\mathrm{Ind}({\uparrow}(x\vee
x^{*}))=\mathrm{Ind}({\uparrow}(x\vee
0_{L_{1}}))=\mathrm{Ind}({\uparrow}x).\] Thus, it suffices to prove
that $\mathrm{Ind}({\uparrow}x)\leqslant k-1$. By Proposition
\ref{pseudo} and our assumption we have that
\[(x,y)\vee (x,y)^{*}=(x,y)\vee(x^{*},y^{*})=(x,y)\vee
(0_{L_{1}},0_{L_{2}})\] and thus, by Lemma \ref{lem1} and Theorem
\ref{mtheo2},
\[\max\{\mathrm{Ind}({\uparrow}x),\mathrm{Ind}({\uparrow}y)\}=\mathrm{Ind}({\uparrow}(x,y))=\mathrm{Ind}({\uparrow}((x,y)\vee
(x,y)^{*}))\] and therefore,
\[\max\{\mathrm{Ind}({\uparrow}x),\mathrm{Ind}({\uparrow}y))\}\leqslant
k-1,\] that is $\mathrm{Ind}({\uparrow}x)\leqslant k-1$. Therefore,
$\mathrm{Ind}(L_1)\leqslant k$.

Similarly, we can prove that $\mathrm{Ind}(L_2)\leqslant k$. Hence,
\[\max\{\mathrm{Ind}(L_1),\mathrm{Ind}(L_2)\}\leqslant k\] and thus, by
Theorem \ref{mtheo2} we have that $\mathrm{Ind}(L_1\times
L_2)\leqslant k$. $\Box$

\begin{corollary}Let $(L_1,\leqslant_1)$ and $(L_2,\leqslant_2)$ be two finite
lattices such that:
\begin{enumerate}
\item the SP-property holds for all ${\uparrow}\kappa$ and
${\uparrow}\lambda$, where $\kappa\in L_1$ and $\lambda\in L_2$ and
\item the pseudocomplement of every $w\in
L_i\setminus\{0_{L_{i}}\}$ is $0_{L_{i}}$, for $i=1,2$.
\end{enumerate}
Then $${\rm Ind}(L_1\times L_2)\leqslant {\rm Ind}(L_1\Box
L_2)\leqslant {\rm Ind}(L_1\times L_2)+1.$$
\end{corollary}
{\bf Proof.} By Theorems \ref{mtheo2} and \ref{mtheo3}, we have
$\mathrm{Ind}(L_1\Box L_2)\leqslant \mathrm{Ind}(L_1\times L_2)+1$.
Moreover, by Theorem \ref{thm10} we have $\mathrm{Ind}(L_1\times
L_2)\leqslant \mathrm{Ind}(L_1\Box L_2)$. $\Box$

\begin{remark}{\rm The dimension Ind of the
Cartesian product $L_1\times L_2$ of two finite lattices $L_1$ and
$L_2$ differs from the dimension Ind of their rectangular product.

\medskip
\noindent (1) We consider the finite lattices, which are represented
in Figures \ref{fig:13} and \ref{fig:14}. We have that
$\mathrm{Ind}(L_1)=\mathrm{Ind}(L_2)=0$, $\mathrm{Ind}(L_1 \Box
L_2)=1$ and by Theorem \ref{mtheo2}, $\mathrm{Ind}(L_1\times
L_2)=0$.

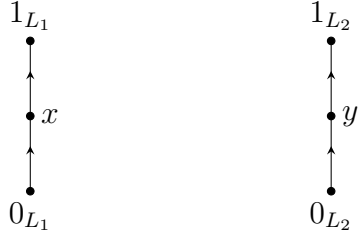
\begin{figure}[H]
\tikzstyle arrowstyle=[scale=1] \tikzstyle
directed=[postaction={decorate,decoration={markings,mark=at position
.6 with {\arrow[arrowstyle]{stealth}}}}]
\begin{center}
\begin{tikzpicture}
\draw[directed] (0,0)--(0,1); \draw[directed] (0,1)--(0,2);
\draw[directed] (4,0)--(4,1); \draw[directed] (4,1)--(4,2);
\draw[fill] (0,0) circle [radius=0.05] node [below] {$0_{L_1}$};
\draw[fill] (0,1) circle [radius=0.05] node [right] {$x$};
\draw[fill] (0,2) circle [radius=0.05] node [above] {$1_{L_1}$};
\draw[fill] (4,0) circle [radius=0.05] node [below] {$0_{L_2}$};
\draw[fill] (4,1) circle [radius=0.05] node [right] {$y$};
\draw[fill] (4,2) circle [radius=0.05] node [above] {$1_{L_2}$};
\end{tikzpicture}
\end{center}
\caption{The finite lattices $(L_1,\leqslant _1)$ and
$(L_2,\leqslant _2)$} \label{fig:13}
\end{figure}
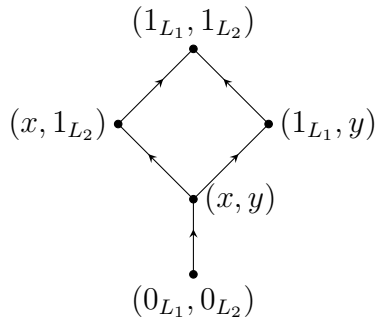
\begin{figure}[H]
\tikzstyle arrowstyle=[scale=1] \tikzstyle
directed=[postaction={decorate,decoration={markings,mark=at position
.6 with {\arrow[arrowstyle]{stealth}}}}]
\begin{center}
\begin{tikzpicture}
\draw[directed] (0,0)--(0,1); \draw[directed] (0,1)--(-1,2);
\draw[directed] (0,1)--(1,2); \draw[directed] (-1,2)--(0,3);
\draw[directed] (1,2)--(0,3); \draw[fill] (0,0) circle [radius=0.05]
node [below] {$(0_{L_1},0_{L_2})$}; \draw[fill] (0,1) circle
[radius=0.05] node [right] {$(x,y)$}; \draw[fill] (-1,2) circle
[radius=0.05] node [left] {$(x,1_{L_2})$}; \draw[fill] (1,2) circle
[radius=0.05] node [right] {$(1_{L_1},y)$}; \draw[fill] (0,3) circle
[radius=0.05] node [above] {$(1_{L_1},1_{L_2})$};
\end{tikzpicture}
\end{center}
\caption{The finite lattice $(L_1 \Box L_2,\leqslant)$}
\label{fig:14}
\end{figure}

\medskip
\noindent (2) We consider the finite lattices of Figure \ref{fig:1}.
Then $\mathrm{Ind}(L_1)=0$ and $\mathrm{Ind}(L_2)=1$ and thus, by
Theorem \ref{mtheo2} we have that $\mathrm{Ind}(L_2\times L_1)=1$.
For the rectangular product $L_2\Box L_1$, which is represented by
the diagram of Figure \ref{fig:15}, we have $\mathrm{Ind}(L_2\Box
L_1)=0$.

Especially, for ${\rm Ind}(L_2\Box L_1)=0$, we follow the argument
(2) of Definition \ref{Ind}. We observe that for every $(a,b)\in
L_2\Box L_1$ and for every $(u.v)\in L_2\Box L_1$ with $(a,b)\vee
(u,v)=(1_{L_{2}},1_{L_{1}})$ we always find an element $(x,y)\in
L_2\Box L_1$ such that $(x,y)\leqslant (u,v)$, $(a,b)\vee
(u,v)=(1_{L_{2}},1_{L_{1}})$ and $${\rm
Ind}({\uparrow}((x,y)^{*}\vee(x,y)))={\rm
Ind}(\{(1_{L_2},1_{L_1})\})=-1.$$ Thus, ${\rm Ind}(L_2\Box L_1)=0$.

\begin{figure}[H]
\tikzstyle arrowstyle=[scale=1] \tikzstyle
directed=[postaction={decorate,decoration={markings,mark=at position
.6 with {\arrow[arrowstyle]{stealth}}}}]
\begin{center}
\begin{tikzpicture}
\draw [directed] (0,1) -- (-3,2); \draw [directed] (0,1) -- (3,2);
\draw [directed] (-3,2) -- (-4,4); \draw [directed] (3,2) -- (4,4);
\draw [directed] (-3,2) -- (-2,4); \draw [directed] (3,2) -- (2,4);
\draw [directed] (2,4) -- (3,6); \draw [directed] (4,4) -- (3,6);
\draw [directed] (-2,4) -- (-3,6); \draw [directed] (-4,4) --
(-3,6); \draw [directed] (0,4) -- (-1,6); \draw [directed] (0,4) --
(1,6); \draw [directed] (-4,4) -- (-1,6); \draw [directed] (-2,4) --
(1,6); \draw [directed] (4,4) -- (1,6); \draw [directed] (2,4) --
(-1,6); \draw [directed] (-3,6) -- (0,7); \draw [directed] (3,6) --
(0,7); \draw [directed] (-1,6) -- (0,7); \draw [directed] (1,6) --
(0,7); \draw [directed] (-3,2) -- (0,4); \draw [directed] (3,2) --
(0,4); \draw[fill] (0,1) circle [radius=0.05] node [below]
{$(0_{L_2},0_{L_1})$}; \draw[fill] (0,4) circle [radius=0.05] node
[left]  {$(y_1,1_{L_1})$}; \draw[fill] (-2,4) circle [radius=0.05]
node [left]  {$(y_3,x_1)$}; \draw[fill] (-3,2) circle [radius=0.05]
node [left]  {$(y_1,x_1)$}; \draw[fill] (3,6) circle [radius=0.05]
node [right]  {$(1_{L_2},x_2)$}; \draw[fill] (-3,6) circle
[radius=0.05] node [left]  {$(1_{L_2},x_1)$}; \draw[fill] (1,6)
circle [radius=0.05] node [right]  {$(y_3,1_{L_1})$}; \draw[fill]
(-4,4) circle [radius=0.05] node [left] {$(y_2,x_1)$}; \draw[fill]
(4,4) circle [radius=0.05] node [right] {$(y_3,x_2)$}; \draw[fill]
(-1,6) circle [radius=0.05] node [left] {$(y_2,1_{L_1})$};
\draw[fill] (3,2) circle [radius=0.05] node [right] {$(y_1,x_2)$};
\draw[fill] (2,4) circle [radius=0.05] node [right] {$(y_2,x_2)$};
\draw[fill] (0,7) circle [radius=0.05] node [above]
{$(1_{L_2},1_{L_1})$};
\end{tikzpicture}
\end{center}
\caption{The finite lattice $(L_2 \Box L_1,\leqslant)$}
\label{fig:15}
\end{figure}
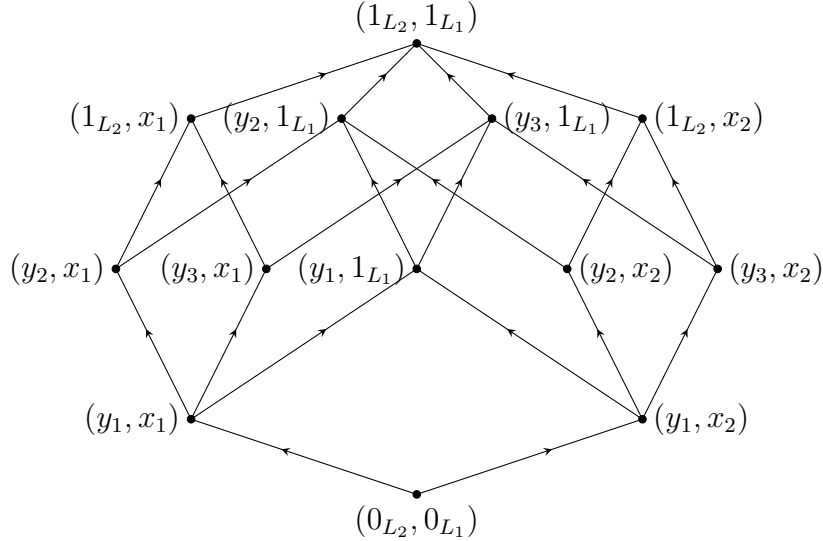}
\end{remark}

\begin{remark}{\rm The dimension Ind of the lexicographic product $L_1 \diamond L_2$ of two
finite lattices $L_1$ and $L_2$ differs from the dimension Ind of
their rectangular product.

\medskip
\noindent (1) We consider the finite lattices $(L_1,\leqslant_1)$
and $(L_3,\leqslant_3)$ of Figures \ref{fig:1} and \ref{fig:9},
respectively. Then as we have seen in Figure \ref{fig:11}, we have
$\mathrm{Ind}(L_3\diamond L_1)=1$. For their rectangular product,
the figure of which is isomorphic to the lattice $(L_1,\leqslant_1)$
of Figure \ref{fig:1}, we have ${\rm Ind}(L_3\Box L_1)=0$.

\medskip
\noindent (2) We consider the lattices $(L_2,\leqslant_2)$ and
$(L_3,\leqslant_3)$ of Figures \ref{fig:1} and \ref{fig:9},
respectively. As we have seen in Figure \ref{fig:12}, we have
$\mathrm{Ind}(L_2\diamond L_3)=0$. For their rectangular product,
the figure of which is isomorphic to the lattice $(L_2,\leqslant_2)$
of Figure \ref{fig:1}, we have ${\rm Ind}(L_2\Box L_3)=1$.}
\end{remark}

\section*{Results and Discussion}

\noindent In this paper, we insert a new dimension, called large
inductive dimension, in the realm of finite lattices. We investigate
several of its properties; sublattice, linear sum, Cartesian,
lexicographic, rectangular product properties. In addition, we study
its relation with known lattice-notions; the small inductive
dimension, the Krull dimension, the covering dimension and the
height, proving the independence of this new dimension in the
chapter of Lattice Dimension Theory.

\section*{Declaration of interests}

\noindent All authors have nothing to declare.

\section*{Declaration of generative AI in scientific writing}

\noindent All authors have nothing to declare.

\section*{Funding sources}

\noindent The authors declare that no funds, grants, or other
support were received during the preparation of this manuscript.

\section*{References}

\bibliographystyle{elsarticle-num-names}

\end{document}